\documentclass{amsart}
\usepackage[toc,page]{appendix}
\usepackage{geometry}
\usepackage[all,cmtip]{xy}
\usepackage{amsmath}
\usepackage{amsfonts}
\usepackage{amssymb}
\usepackage{mathrsfs}   
\usepackage{amsthm}
\usepackage{xcolor}
\usepackage{cite}
\usepackage{stmaryrd}
\usepackage{graphicx}
\usepackage{pgf}
\usepackage{eepic}
\usepackage{pstricks}
\usepackage[all,cmtip]{xy}
\usepackage{epsfig}
\usepackage{fancyhdr}
\usepackage[backref=page]{hyperref}
\renewcommand*{\backref}[1]{}
\renewcommand*{\backrefalt}[4]{({%
    \ifcase #1 Not cited.%
          \or page~#2%
          \else pages #2%
    \fi%
    })}
    
 \usepackage{mathptmx}  

 \newcommand\imCMsym[4][\mathord]{%
  \DeclareFontFamily{U} {#2}{}
  \DeclareFontShape{U}{#2}{m}{n}{
    <-6> #25
    <6-7> #26
    <7-8> #27
    <8-9> #28
    <9-10> #29
    <10-12> #210
    <12-> #212}{}
  \DeclareSymbolFont{CM#2} {U} {#2}{m}{n}
  \DeclareMathSymbol{#4}{#1}{CM#2}{#3}
}
\newcommand\alsoimCMsym[4][\mathord]{\DeclareMathSymbol{#4}{#1}{CM#2}{#3}}

\imCMsym{cmmi}{124}{\CMjmath}
\imCMsym[\mathop]{cmsy}{113}{\CMamalg}
\imCMsym[\mathop]{cmex}{96}{\CMcoprod}
\alsoimCMsym[\mathop]{cmex}{97}{\CMbigcoprod}
 
\theoremstyle{plain}
\newtheorem*{theoremu}{Theorem}
\newtheorem{theorem}{Theorem}[section]

\newtheorem{proposition}[theorem]{Proposition}
\newtheorem{corollary}[theorem]{Corollary}

\newtheorem{lemma}[theorem]{Lemma}

\theoremstyle{definition}

\newtheorem{hypothesis}[theorem]{Hypothesis}
\newtheorem{definition}[theorem]{Definition}

\theoremstyle{remark}
\newtheorem{remark}[theorem]{Remark}
\newtheorem{example}[theorem]{Example}

\newcommand{\N}{{\mathbb N}}
\newcommand{\Z}{{\mathbb Z}}
\newcommand{\Q}{{\mathbb Q}}

\newcommand{\C}{{\mathbb C}}

\newcommand{\F}{{\mathbb F}}

\newcommand{\eps}{\varepsilon}%

\newcommand{\nic}[1]{\cur{N}\mathrm{IC}(#1)}
\newcommand{\spec}[1]{\mathrm{Spec}\left(#1\right)}
\newcommand{\cur}[1]{\mathcal{#1}}

\newcommand{\isomto}{\overset{\sim}{\rightarrow}}

\newcommand{\RdR}[1]{\mathbf{R}^{#1}_\mathrm{dR}f_*}

\newcommand{\rig}{\mathrm{rig}}

\title{Relative fundamental groups and rational points}

\author{Christopher Lazda}
       \address{Dipartimento di Matematica Pura e Applicata \\
        Torre Archimede, Via Trieste, 63 \\ 
        35121 Padova \\ 
        Italia}
       \email{lazda@math.unipd.it}

\newcommand{\lie}[1]{\mathfrak{#1}}

\begin{document}

\begin{abstract} In this paper we define a relative rigid fundamental group, which associates to a section $p$ of a smooth and proper morphism $f:X\rightarrow S$ in characteristic $p$, with $\dim S=1$, a Hopf algebra in the ind-category of overconvergent $F$-isocrystals on $S$. We prove a base change property, which says that the fibres of this object are the Hopf algebras of the rigid fundamental groups of the fibres of $f$. We explain how to use this theory to define period maps as Kim does for varieties over number fields, and show in certain cases that the targets of these maps can be interpreted as varieties.
\end{abstract}

\maketitle 

\tableofcontents

\section*{Introduction}

Let $K$ be a number field and let $C/K$ be a smooth, projective curve of genus $g>1$, with Jacobian $J$. Then a famous theorem of Faltings states that the set $C\left(K\right)$ of $K$-rational points on $C$ is finite. The group $J\left(K\right)$ is finitely generated, and under the assumption that its rank is strictly less than $g$, Chabauty in \cite{Cha41} was able to prove this theorem using elementary methods as follows. Let $v$ be a place of $K$, of good reduction for $C$, and denote by $C_v,J_v$ the base change to $K_v$. Then Chabauty defines a homomorphism
\begin{equation}\log : J\left(K_v\right)\rightarrow H^0\left(J_v,\Omega^1_{J_v/K_v}\right)
\end{equation}
and shows that there exists a non zero linear functional on $H^0(J_v,\Omega^1_{J_v/K_v})$ which vanishes on the image of $J\left(K\right)$. He then proves that pulling this back to $J\left(K_v\right)$ gives an analytic function on $J\left(K_v\right)$, which is not identically zero on $C\left(K_v\right)$, and which vanishes on $J\left(K\right)$. Hence $C\left(K\right)\subset C\left(K_v\right)\cap J\left(K\right)$ must be finite as it is contained in the zero set of a non-zero analytic function on $C\left(K_v\right)$.

In \cite{Kim09}, Kim describes what he calls a `non-abelian lift' of this method. Fix a point $p\in C\left(K\right)$. By considering the Tannakian category of integrable connections on $C_v$, one can define a `de Rham fundamental group' $U^\mathrm{dR}=\pi_1^\mathrm{dR}\left(C_v,p\right)$, which is a pro-unipotent group scheme over $K_v$, as well as, for any other $x\in C\left(K_v\right)$, path torsors $P^\mathrm{dR}\left(x\right)=\pi_1^\mathrm{dR}\left(C_v,x,p\right)$ which are right torsors under $U^\mathrm{dR}$. These group schemes and torsors come with extra structure, namely that of a Hodge filtration and, by comparison with the crystalline fundamental group of the reduction of $C_v$, a Frobenius action. He then shows that such torsors are classified by $U^\mathrm{dR}/F^0$, and hence one can define `period maps'
\begin{equation} j_n:C\left(K_v\right)\rightarrow U^\mathrm{dR}_n/F^0
\end{equation}
where $U_n^\mathrm{dR}$ is the $n$th level nilpotent quotient of $U^\mathrm{dR}$. If $n=2$ then $j_n$ is just the composition of the above $\log$ map with the inclusion $C\left(K_v\right)\rightarrow J\left(K_v\right)$. By analysing the image of this map, he is able to prove finiteness of $C\left(K\right)$ under certain conditions, namely if the dimension of $U_n^\mathrm{dR}/F^0$ is greater than the dimension of the target of a global period map defined using the category of lisse \'{e}tale sheaves on $C$. Moreover, when $n=2$, this condition on dimensions is essentially  Chabauty's condition that $\mathrm{rank}_\Z J\left(K\right)<\mathrm{genus}\left(C\right)$ (modulo the Tate-Shafarevich conjecture).

Our interest lies in trying to develop a function field analogue of these ideas. The analogy between function fields in one variable over finite fields and number fields has been a fruitful one throughout modern number theory, and indeed the analogue of Mordell's conjecture was first proven for function fields by Grauert. In this paper we discuss the problem of defining a good analogue of the global period map. This is defined in \cite{Kim09} using the Tannakian category of lisse $\Q_p$ sheaves on $X$, and this approach will not work in the function field setting. Neither $p$-adic nor $\ell$-adic \'{e}tale cohomology will give satisfactory answers, the first because, for example, the resulting fundamental group will be moduli dependent, i.e. will not be locally constant in families (see for example \cite{Tam04}), and the second because the $\ell$-adic topology on the resulting target spaces for period maps will not be compatible with the $p$-adic topology on the source varieties. Instead we will work with the category of overconvergent $F$-isocrystals.

Let $K$ be a finite extension of $\F_p\left(t\right)$, and let $k$ be the field of constants of $K$, i.e. the algebraic closure of $\F_p$ inside $K$. Let $\overline{S}$ be the unique smooth projective, geometrically irreducible curve over $k$ whose function field is $K$. If $C/K$ is a smooth, projective, geometrically integral curve then one can choose a regular model for $C$. This is a regular, proper surface $\overline{X}/k$, equipped with a flat, proper morphism $f:\overline{X}\rightarrow \overline{S}$ whose generic fibre is $C/K$. Let $S\subset \overline{S}$ be the smooth locus of $f$, and denote by $f$ also the pullback $f:X\rightarrow S$. The idea is to construct, for any section $p$ of $f$, a `non-abelian isocrystal' on $S$ whose fibre at any closed point $s$ `is' the rigid fundamental group $\pi_1^\mathrm{rig}\left(X_s,p_s\right)$. The idea behind how to construct such an object is very simple.

Suppose that $f:X\rightarrow S$ is a Serre fibration of topological spaces, with connected base and fibres. If $p$ is a section, then for any $s\in S$ the homomorphism $\pi_1\left(X,p\left(s\right)\right)\rightarrow \pi_1\left(S,s\right)$ is surjective, and $\pi_1\left(S,s\right)$ acts on the kernel via conjugation. This corresponds to a locally constant sheaf of groups on $S$, and the fibre over any point $s\in S$ is just the fundamental group of the fibre $X_s$. This approach makes sense for any fundamental group defined algebraically as the Tannaka dual of a category of `locally constant' coefficients. So if $f:X\rightarrow S$ is a morphism of smooth varieties with section $p$, then $f_*:\pi_1^{\cur{C}_X}\left(X,x\right)\rightarrow \pi_1^{\cur{C}_S}\left(S,s\right)$ is surjective, and $\pi_1^{\cur{C}_S}\left(S,s\right)$ acts on the kernel. Here $\cur{C}_{\left(-\right)}$ is any appropriate category of coefficients, for example vector bundles with integrable connection, unipotent isocrystals etc., and e$.$g$.$ $\pi_1^{\cur{C}_X}\left(X,x\right)$ is the Tannaka dual of this category with respect to the fibre functor $x^*$. This gives the kernel of $f_*$ the structure of an `affine group scheme over $\cur{C}_S$', and it makes sense to ask what the fibre is over any closed point $s\in S$. The main theorem of the first chapter is the following.

\begin{theoremu} Suppose that $f:X\rightarrow S$ is a smooth morphism of smooth varieties over an algebraically closed field $k$ of characteristic zero. Assume that $f$ has geometrically connected fibres, and that $S$ is a geometrically connected affine curve. Assume further that $X$ is the complement of a relative normal crossings divisor in a smooth and proper $S$-scheme $\overline{X}$. Let $\cur{C}_S$ be the category of vector bundles with a regular integrable connection on $S$, and let $\cur{C}_X$ be the category of vector bundles with a regular integrable connection on $X$ which are iterated extensions of those of the form $f^*\cur{E}$, with $\cur{E}\in \cur{C}_S$. Then the fibre of the corresponding affine group scheme over $\cur{C}_S$ at $s\in S$ is the de Rham fundamental group $\pi_1^\mathrm{dR}\left(X_s,p_s\right)$ of the fibre.
\end{theoremu}

Thus with strong hypotheses on the base $S$, we have a good working definition of a relative fundamental group. We would ideally like to remove these hypotheses, and it seems as though a good way to do this would be to use the methods of `relative rational homotopy theory' similar to Navarro-Aznar's work in \cite{NA93}. In positive characteristic at least, this approach will be taken up in future work.

In Chapter 2 we discuss path torsors in the relative setting. We show in particular that for any other section $q$ of $f$ one can define an affine scheme $\pi_1^\mathrm{dR}\left(X/S,q,p\right)$ over $\cur{C}_S$ which is a right torsor under the relative de Rham fundamental group $\pi_1^\mathrm{dR}\left(X/S,p\right)$. The upshot of this is that we obtain
\begin{equation} j_n :X\left(S\right)\rightarrow H^1\left(S,\pi_1^\mathrm{dR}\left(X/S,p\right)_n\right)
\end{equation}
which are a coarse characteristic zero function field analogue of Kim's global period maps. Of course, if we were really interested in the characteristic zero picture, we would want to define Hodge structures on these objects, and thus obtain finer period maps. However, our main interest lies in the positive characteristic case, and so we don't pursue these questions.

 In Chapter 3 we define the relative rigid fundamental group in positive characteristic, mimicking the definition in characteristic zero. Instead of the category of vector bundles with regular integrable connections, we consider the category of overconvergent $F$-isocrystals (throughout Chapter 3 we will be over a finite field, and Frobenius will always mean the \emph{linear} Frobenius). We then proceed to use Caro's theory of cohomological operations for arithmetic $\cur{D}$-modules in order to prove the analogue of the above theorem in positive characteristic. Although sufficient for our ultimate end goal, where our bases are geometrically connected, smooth curves, it would be pleasing to have a formalism that worked in greater generality. As mentioned above, this will form part of a future work.

The upshot of this is that for a smooth and proper map $f:X\rightarrow S$ with geometrically connected fibres, $S$ a smooth, geometrically connected curve over a finite field $k$, and a section $p$ of $f$, we can define an affine group scheme $\pi_1^\rig(X/S,p)$ over the category of overconvergent $F$-isocrystals on $S$, which we call the relative fundamental group at $p$. The fibre of this over any point $s\in S$ is just the unipotent rigid fundamental group of the fibre $X_s$ of $f$ over $s$. As in the zero characteristic case, the general Tannakian formalism gives us path torsos $\pi_1^\rig(X/S,p,q)$ for any other $q\in X(S)$, and hence we can define a period map
\begin{equation} X(S)\rightarrow H^1_{F,\rig}(S,\pi_1^\rig(X/S,p))
\end{equation} 
where the RHS is a classifying set of $F$-torsors under $\pi_1^\rig(X/S,p)$, as well as finite level versions given by pushing out along the quotient map $\pi_1^\rig(X/S,p)\rightarrow\pi_1^\rig(X/S,p)_n$.

Finally, we study the targets of these period maps, and show that after replacing $H^1_{F,\rig}(S,\pi_1^\rig(X/S,p))$, the set classifying $F$-torsors, by $H^1_{\rig}(S,\pi_1^\rig(X/S,p))^{\phi=\mathrm{id}}$, the Frobenius invariant part of the set classifying torsors without $F$-structure, then under very restrictive hypotheses on the morphism $f:X\rightarrow S$, we obtain the structure of an algebraic variety. The argument here is just a translation of the original argument of Kim into our context, and what for us are restrictive hypotheses are automatically satisfied in his case.

We are still a long way away from getting a version of Kim's methods to work for function fields. There is still the question of how to define the analogue of the local period maps, and also to show that the domains of the period maps have the structure of varieties. Even then, it is very unclear what the correct analogue of the local integration theory will be in positive characteristic. There is still a very large amount of work to be done if such a project is to be completed.

\section{Relative de Rham fundamental groups}\label{zero}

Let $f:X\rightarrow S$ be a smooth morphism of smooth complex varieties, and suppose that $f$ admits a good compactification, that is, there exists $\overline{X}$ smooth and proper over $S$, an open immersion $X\hookrightarrow \overline{X}$ over $S$, such that $D=\overline{X}\setminus X$ is a relative normal crossings divisor in $\overline{X}$. Let $p\in X\left(S\right)$ be a section. For every closed point $s\in S$ with fibre $X_s$, one can consider the topological fundamental group $G_s:=\pi_1\left(X_s^\mathrm{an},p\left(s\right)\right)$, and as $s$ varies, these fit together to give a locally constant sheaf $\pi_1\left(X/S,p\right)$ on $S^\mathrm{an}$. Let
\begin{equation}
\hat{\cur{U}}\left(\mathrm{Lie}\; G_s\right):=\varprojlim \C[G_s]/\lie{a}^n
\end{equation}  
denote the completed enveloping algebra of the Malcev Lie algebra of $G_s$, where $\lie{a}\subset \C[G_s]$ is the augmentation ideal. According to Proposition 4.2 of \cite{HZ87}, as $s$ varies, these fit together to give a pro-local system on $S^\mathrm{an}$, i.e. a pro-object $\hat{\cur{U}}_p^\mathrm{top}$ in the category of locally constant sheaves of finite dimensional $\C$-vector spaces on $S^\mathrm{an}$. (Their theorem is a lot stronger than this, but this is all we need for now). According to Th\'{e}or\`{e}me 5.9 in Chapter II of \cite{Del70}, the pro-vector bundle with integrable connection $\hat{\cur{U}}_p^\mathrm{top}\otimes_\C \cur{O}_{S^\mathrm{an}}$ has a canonical algebraic structure. Thus given a smooth morphism $f:X\rightarrow S$ as above, with section $p$, one can construct a pro-vector bundle with connection $\hat{\cur{U}}_p$ on $S$, whose fibre at any closed point $s\in S$ is the completed enveloping algebra of the Malcev Lie algebra of $\pi_1\left(X_s^\mathrm{an},p\left(s\right)\right)$.

Denoting by $\lie{g}_s$ the Malcev Lie algebra of $\pi_1\left(X_s^\mathrm{an},p\left(s\right)\right)$, $\hat{\cur{U}}\left(\lie{g}_s\right)=(\hat{\cur{U}}_p)_s$ can be constructed algebraically, as $\lie{g}_s$ is equal to $\mathrm{Lie}\;\pi_1^\mathrm{dR}\left(X_s,p_s\right)$, the Lie algebra of the Tannaka dual of the category of unipotent vector bundles with integrable connection on $X_s$. This suggests the question of whether or not there is an algebraic construction of $\hat{\cur{U}}_p$?

We will not directly answer this question - instead we will construct the Lie algebra associated to $\hat{\cur{U}}_p$ - this is a pro-system $\hat{\cur{L}}_p$ of Lie algebras with connection on $S$. The way we will do so is very simple, and is closely related to ideas used in \cite{Wil97} to study relatively unipotent mixed motivic sheaves.

\begin{definition} To save ourselves saying the same thing over and over again, we make the following definition. A `good' morphism is a smooth morphism $f:X\rightarrow S$ of smooth varieties over a field $k$, with geometrically connected fibres and base, such that $X$ is the complement of a relative normal crossings divisor in a smooth, proper $S$-scheme $\overline{X}$. Throughout this section we will assume that the ground field $k$ is algebraically closed of characteristic $0$.
\end{definition}

We will assume that the reader is familiar with Tannakian categories, a good introductory reference is \cite{MD81}. If $\cur{T}$ is a Tannakian category over a field $k$, and $\omega$ is a fibre functor on $\cur{T}$, in the sense of \S1.9 of \cite{Del90}, we will denote the group scheme representing tensor automorphisms of $\omega$ by $G(\cur{T},\omega)$. We will also use the rudiments of algebraic geometry in Tannakian categories, as explained in \S5 of \cite{Del89} - in particular we will talk about affine (group) schemes over Tannakian categories. We will denote the fundamental groupoid of a Tannakian category by $\pi(\cur{T})$, this is an affine group scheme over $\cur{T}$ which satisfies $\omega(\pi(\cur{T}))=G(\cur{T},\omega)$ for every fibre functor $\omega$ (see for example 6.1 of \cite{Del89}). If $\cur{T}$ is a Tannakian category over $k$, and $k'/k$ is a finite extension, then we will denote the category of $k'$-modules in $\cur{T}$ by either $\cur{T}\otimes_k k'$, or $\cur{T}_{k'}$.

We will also assume familiarity with the theory of integrable connections and regular holonomic $\mathcal{D}$-modules on $k$-varieties, and will generally refer to \cite{Del70} and \cite{Bor87} for details. We say that a regular integrable connection on $X$ is unipotent if it is a successive extension of the trivial connection, and these form a Tannakian subcategory $\cur{N}\mathrm{IC}(X)\subset \mathrm{IC}(X)$ of the Tannakian category of regular integrable connections.

\begin{definition} For $X/k$ smooth and connected, the algebraic and de Rham fundamental groups of $X$ at a closed point $x\in X$ are defined by
\begin{align}\pi_1^\mathrm{alg}\left(X,x\right)&:=x^*\left(\pi\left(\mathrm{IC}\left(X\right)\right)\right)=G\left(\mathrm{IC}\left(X\right),x^*\right)\\
 \pi_1^\mathrm{dR}\left(X,x\right)&:=x^*\left(\pi\left(\cur{N}\mathrm{IC}\left(X\right)\right)\right) =G\left(\cur{N}\mathrm{IC}\left(X\right),x^*\right). \end{align} 
\end{definition}

\begin{remark} \label{comp1}It follows from the Riemann-Hilbert correspondence that if $k=\C$, then these affine group schemes are the pro-algebraic and pro-unipotent completions of $\pi_1\left(X^\mathrm{an},x\right)$ respectively.
\end{remark}

If $f:X\rightarrow Y$ is a morphism of smooth $k$-varieties, then we can form the pullback of vector bundles with integrable connection on $Y$, which preserves regularity and is the usual pull-back on the underlying $\cur{O}_Y$-module. This induces a homomorphism $f_*:\pi_1^\#(X,x)\rightarrow \pi_1^\#(Y,f(x))$ for $\#=\mathrm{dR},\mathrm{alg}$.

\subsection{The relative fundamental group and its pro-nilpotent Lie algebra}\label{LIE}

Let $f:X\rightarrow S$ be a `good' morphism. A regular integrable connection $E$ on $X$ is said to be relatively unipotent if there exists a filtration by horizontal sub-bundles, whose graded objects are all in the essential image of $f^*:\mathrm{IC}\left(S\right)\rightarrow \mathrm{IC}\left(X\right)$. We will denote the full subcategory of relatively unipotent objects in $\mathrm{IC}\left(X\right)$ by $\cur{N}_f\mathrm{IC}\left(X\right)$, which is a Tannakian subcategory. Suppose that $p\in X\left(S\right)$ is a section of $f$. We have functors of Tannakian categories
\begin{equation}\xymatrix{
\cur{N}_f\mathrm{IC}\left(X\right)\ar@<0.6ex>[r]^-{p^*}&\ar@<0.6ex>[l]^-{f^*} \mathrm{IC}\left(S\right)
}\end{equation}  
and hence, after choosing a point $s\in S(k)$, homomorphisms
\begin{equation}\xymatrix{
G\left(\cur{N}_f\mathrm{IC}\left(X\right),p\left(s\right)^*\right)\ar@<0.6ex>[r]^-{f_*}&\ar@<0.6ex>[l]^-{p_*} G\left(\mathrm{IC}\left(S\right),s^*\right)
}\end{equation}  
between their Tannaka duals. Let $K_s$ denote the kernel of $f_*$. Then the splitting $p_*$ induces an action of $\pi_1^\mathrm{alg}\left(S,s\right)=G\left(\mathrm{IC}\left(S\right),s^*\right)$ on $K_s$ via conjugation. This corresponds to an affine group scheme over $\mathrm{IC}\left(S\right)$.

\begin{lemma} This affine group scheme is independent of $s$.
\end{lemma}

\begin{proof} Thanks to \cite{Del89}, \S6.10, $f_*, p_*$ above come from homomorphisms 
\begin{equation}\xymatrix{
p^*\left(\pi\left(\cur{N}_f\mathrm{IC}\left(X\right)\right)\right)\ar@<0.6ex>[r]^-{f_*}&\ar@<0.6ex>[l]^-{p_*} \pi\left(\mathrm{IC}\left(S\right)\right)
}\end{equation}  
of affine group schemes over $\mathrm{IC}\left(S\right)$. If we let $\cur{K}$ denote the kernel of $f_*$, then $K_s=s^*\left(\cur{K}\right)$.
\end{proof}

\begin{definition} The relative de Rham fundamental group $\pi_1^\mathrm{dR}\left(X/S,p\right)$ of $X/S$ at $p$ is defined to be the affine group scheme $\cur{K}$ over $\mathrm{IC}(S)$.
\end{definition}

Let $i_s:X_s\rightarrow X$ denote the inclusion of the fibre over $s$. Then there is a canonical functor $i_s^*:\cur{N}_f\mathrm{IC}\left(X\right)\rightarrow \cur{N}\mathrm{IC}\left(X_s\right)$. This induces a homomorphism $\pi_1^\mathrm{dR}\left(X_s,p_s\right)\rightarrow G\left(\cur{N}_f\mathrm{IC}\left(X\right),p_s^*\right)$ which is easily seen to factor through the fibre $\pi_1^\mathrm{dR}\left(X/S,p\right)_s:=s^*(\cur{K})=K_s$ of $\pi_1^\mathrm{dR}\left(X/S,p\right)$ over $s$.

\begin{theorem} \label{woot} Suppose that $k=\C$. Then $\phi:\pi_1^\mathrm{dR}\left(X_s,p_s\right)\rightarrow \pi_1^\mathrm{dR}\left(X/S,p\right)_s$ is an isomorphism.
\end{theorem}

\begin{proof} The point $s$ gives us fibre functors $p_s^*$ on $\nic{X_s}$, $p\left(s\right)^*$ on $\cur{N}_f\mathrm{IC}\left(X\right)$ and $s^*$ on $\mathrm{IC}\left(S\right)$. Write
\begin{equation}
\cur{K} = G\left(\nic{X_s},p_s^*\right), \quad
\cur{G} = G\left(\cur{N}_f\mathrm{IC}\left(X\right),p\left(s\right)^*\right), \quad
\cur{H} = G\left(\mathrm{IC}\left(S\right),s^*\right)
\end{equation}
and also let 
\begin{equation}
K=\pi_1\left(X_s^\mathrm{an},p\left(s\right)\right), \quad
G= \pi_1\left(X^\mathrm{an},p\left(s\right)\right),\quad
H=\pi_1\left(S^\mathrm{an},s\right)
\end{equation}
be the topological fundamental groups of $X_s,X,S$ respectively. Then $\cur{K}=K^\mathrm{un}$, the pro-unipotent completion of $K$, and $\cur{H}=H^\mathrm{alg}$, the pro-algebraic completion of $H$. We need to show that the sequence of affine group schemes 
\begin{equation} 1\rightarrow \cur{K}\rightarrow \cur{G}\rightarrow \cur{H}\rightarrow 1
\end{equation}  
is exact, and we will use the equivalences of categories
\begin{equation} \mathrm{IC}\left(X\right)\isomto\mathrm{Rep}_{\C}\left(\pi_1\left(X^\mathrm{an},p\left(s\right)\right)\right),\quad  
\mathrm{IC}\left(S\right)\isomto\mathrm{Rep}_{\C}\left(\pi_1\left(S^\mathrm{an},s\right)\right)
\end{equation}
\begin{equation}
\mathrm{IC}\left(X_s\right)\isomto\mathrm{Rep}_{\C}\left(\pi_1\left(X_s^\mathrm{an},p\left(s\right)\right)\right). \end{equation}
By Proposition 1.3 in Chapter I of \cite{Wil97}, $\ker\left(\cur{G}\rightarrow \cur{H}\right)$ is pro-unipotent. Hence according to Proposition 1.4 of \emph{loc. cit.}, in order to show that $\phi$ is an isomorphism, we must show the following.

\begin{itemize} 
\item If $E\in\cur{N}_f\mathrm{IC}\left(X\right)$ is such that $i_s^*\left(E\right)$ is trivial, then $E\cong f^*\left(
F\right)$ for some $F$ in $\mathrm{IC}\left(S\right)$.
\item Let $E\in\cur{N}_f\mathrm{IC}(X)$, and let $F_0\subset i_s^*(E)$ denote the largest trivial sub-object. Then there exists $E_0\subset E$ such that $F_0=i_s^*(E_0)$.
\item There is a pro-action of $\cur{G}$ on $\hat{\cur{U}}\left(\mathrm{Lie}\; \cur{K}\right)$ such that the corresponding action of $\mathrm{Lie}\; \cur{G}$ extends the left multiplication by $\mathrm{Lie} \;\cur{K}$. 
\end{itemize}

The first is straightforward. Since $f$ is topologically a fibration with section $p$, we have a split exact sequence 
\begin{equation} 1\rightarrow K\rightarrow G\leftrightarrows H\rightarrow1 \end{equation}   and a representation $V$ of $G$ such that $K$ acts trivially. We must show that $V$ is the pullback of an $H$-representation - this is obvious! The second is no harder, we must show that if $V$ is a $G$-representation, then $V^K$ is a sub-$G$-module of $V$. But since $K$ is normal in $G$, this is clear. For the third, note that $\hat{\cur{U}}\left(\mathrm{Lie}\;\cur{K}\right)=\hat{\cur{U}}\left(\mathrm{Lie}\;K\right)=\varprojlim \C[K]/\lie{a}^n$, where $\lie{a}$ is the augmentation ideal of $\C[K]$. Let $H$ act on $\C[K]/\lie{a}^n$ by conjugation and $K$ by left multiplication. I claim that $\C[K]/\lie{a}^n$ is finite dimensional, and unipotent as a $K$-representation.

Indeed,  There are extensions of $K$-representations
\begin{equation} 0\rightarrow \lie{a}^{n}/\lie{a}^{n+1}\rightarrow \C[K]/\lie{a}^{n+1}\rightarrow \C[K]/\lie{a}^{n}\rightarrow 0
\end{equation}  
and hence, since the action of $K$ on $\lie{a}^n/\lie{a}^{n+1}$ is trivial, it follows by induction that each $\C[K]/\lie{a}^n$ is unipotent. There are also surjections
\begin{equation} \left(\lie{a}/\lie{a}^2\right)^{\otimes n}\twoheadrightarrow \lie{a}^{n}/\lie{a}^{n+1}
\end{equation}  
for each $n$, and hence by induction, to show finite dimensionality it suffices to show that $\lie{a}/\lie{a}^2$ is finite dimensional. But $\lie{a}/\lie{a}^2\cong K^\mathrm{ab}\otimes_{\Z} \C$ is finite dimensional, as $K$ is finitely generated.

Now, since $\C[K]/\lie{a}^n$ is unipotent as a $K$-representation,it is relatively unipotent as a $G=K\rtimes H$-representation, hence $\C[K]/\lie{a}^n$ is naturally an object in $\mathrm{Rep}_\C\left(\cur{G}\right)$. Thus there is a pro-action of $\cur{G}$ on $\hat{\cur{U}}\left(\mathrm{Lie}\;\cur{K}\right)$, and the action extends left multiplication by $\mathrm{Lie}\;\cur{K}$ as required.
\end{proof}

\begin{remark} The co-ordinate algebra of $\pi_1^\mathrm{dR}(X/S,p)$ is an ind-object in the category of regular integrable connections on $S$. Hence we may view $\pi_1^\mathrm{dR}(X/S,p)$ as an affine group scheme over $S$ in the usual sense, together with a regular integrable connection on the associated $\cur{O}_S$-Hopf algebra.
\end{remark}

If $g:T\rightarrow S$ is any morphism of smooth varieties over $k$, then there is a homomorphism of fundamental groups
\begin{equation}\label{bsemor} \pi_1^\mathrm{dR}\left(X_T/T,p_T\right)\rightarrow \pi_1^\mathrm{dR}\left(X/S,p\right)\times_S T:=g^*(\pi_1^\mathrm{dR}\left(X/S,p\right))
\end{equation}  
which corresponds to a horizontal morphism
\begin{equation} \cur{O}_{\pi_1^\mathrm{dR}\left(X/S,p\right)} \otimes_{\cur{O}_S}\cur{O}_T \rightarrow \cur{O}_{\pi_1^\mathrm{dR}\left(X_T/T,p_T\right)}.
\end{equation}  

\begin{proposition} \label{bse}If $k=\C$ then this is an isomorphism.
\end{proposition}

\begin{proof} We know by the previous theorem that this induces an isomorphism on fibres over any point $t\in T(\C)$. Hence by rigidity, it is an isomorphism.
\end{proof}

Write $G=\pi_1^\mathrm{dR}\left(X/S,p\right)$ and let $G_n$ denote the quotient of $G$ by the $n$th term in its lower central series. Let $A_n$ denote the Hopf algebra of $G_n$, and $I_n\subset A_n$ the augmentation ideal. $L_{n}:=\cur{H}\mathrm{om}_{\cur{O}_S}\left(I_n/I_n^2,\cur{O}_S\right)$ is the Lie algebra of $G_n$. This is a coherent, nilpotent Lie algebra with connection, i.e the bracket $[\cdot,\cdot]:L_{n}\otimes L_{n}\rightarrow L_{n}$ is horizontal. There are natural morphisms $L_{n+1}\rightarrow L_{n}$, which form a pro-system of nilpotent Lie algebras with connection $\hat{L}_p$, whose universal enveloping algebra is the object $\hat{\cur{U}}_p$ considered in the introduction to this section.

\subsection{Towards an algebraic proof of Theorem \ref{woot}}

Although we have a candidate for the relative fundamental group of a `good' morphism $f:X\rightarrow S$ at a section $p$, we have only proved it is a good candidate when the ground field is the complex numbers. One might hope to be able to reduce to the case $k=\C$ via base change and finiteness arguments, but this approach will not work in a straightforward manner. Also, such an argument will not easily adapt to the case of positive characteristic, as in general one will not be able to lift a smooth proper family, even locally on the base. Instead we seek a more algebraic proof. Recall that we have an affine group scheme $\pi_1^\mathrm{dR}\left(X/S,p\right)$ over $\mathrm{IC}(S)$, and a comparison morphism
\begin{equation}\label{comp2}\phi:\pi_1^\mathrm{dR}\left(X_s,p_s\right)\rightarrow \pi_1^\mathrm{dR}\left(X/S,p\right)_s
\end{equation}  
for any point $s\in S$. We want to show that when $S$ is an affine curve, this map is an isomorphism.

It follows from Proposition 1.4 in Chapter I of \cite{Wil97} and Appendix A of \cite{EHS07} that we need to prove the following:
\begin{itemize}
\item (Injectivity) Every $E\in\cur{N}\mathrm{IC}\left(X_s\right)$ is a sub-quotient of $i^*_s\left(F\right)$ for some $F\in\cur{N}_f\mathrm{IC}\left(X\right)$.
\item (Surjectivity I) \label{surj} Suppose that $E\in\cur{N}_f\mathrm{IC}\left(X\right)$ is such that $i^*_s\left(E\right)$ is trivial. Then there exists $F\in\mathrm{IC}(S)$ such that $E\cong f^*\left(F\right)$.
\item (Surjectivity II) \label{surj2} Let $E\in\cur{N}_f\mathrm{IC}(X)$, and let $F_0\subset i_s^*(E)$ denote the largest trivial sub-object. Then there exists $E_0\subset E$ such that $F_0=i_s^*(E_0)$.

\end{itemize}

To do so, we will need to use the language of algebraic $\cur{D}$-modules. We define the functor 
$$ 
f_*^\mathrm{dR}:\cur{N}_f\mathrm{IC}(X)\rightarrow \mathrm{IC}(S)
$$
by $f_*^\mathrm{dR}(E)=\cur{H}^{-d}(f_+E) $ where $f_+$ is the usual push-forward for regular holonomic complexes of $\cur{D}$-modules, $d$ is the relative dimension of $f:X\rightarrow S$, and we are considering a regular integrable connection on $X$ as a $\cur{D}_X$-module in the usual way. 

\begin{lemma} The functor $f_*^\mathrm{dR}$ lands in the category of regular integrable connections, and  is a right adjoint to $f^*$.
\end{lemma}

\begin{proof}
The content of the first claim is in the coherence of direct images in de Rham cohomology, using the comparison result 1.4 of \cite{DMSS00}, and the fact that a regular holonomic $\cur{D}_X$-module is a vector bundle iff it is coherent as an $\cur{O}_X$-module.

To see this coherence, we first use adjointness of $f_+$ and $f^+$, together with the facts that $f^+\cur{O}_S=\cur{O}_X[-d]$ and $f_+\cur{O}_X$ is concentrated in degrees $\geq -d$, to get  canonical adjunction morphism $f_*^\mathrm{dR}(\cur{O}_X)\rightarrow \cur{O}_S$ of regular holonomic $\cur{D}_X$-modules. This is an isomorphism by base changing to $\C$ and comparing with the usual topological push-forward of the constant sheaf $\C$. Hence $f_*^\mathrm{dR}\cur{O}_X$ is coherent, and via the projection formula, so is $f_*^\mathrm{dR}(f^*F)$ for any $F\in\mathrm{IC}(S)$. Hence using exact sequences in cohomology and induction on unipotence degree, $f_*^\mathrm{dR}E$ is coherent whenever $E$ is relatively unipotent.

To prove to the second claim, we just use that $f^+$ is adjoint to $f_+$, $f^+=f^*[-d]$ on the subcategory of regular integrable connections, and $f_+E$ is concentrated in degrees $\geq-d$ whenever $E$ is a regular integrable connection.
\end{proof}

\begin{remark} Although the Proposition is stated in \cite{DMSS00} for $k=\C$, the same proof works for any algebraically closed field of characteristic zero.
\end{remark}

Thus we get a canonical morphism $\varepsilon_E:f^*f_*^\mathrm{dR}E\rightarrow E$ which is the counit of the adjunction between $f^*$ and $f_*^\mathrm{dR}$.

\begin{example}  \label{explad} Suppose that $S=\mathrm{Spec}\left(k\right)$. Then 
\begin{equation}f_*^\mathrm{dR}E=H^0_\mathrm{dR}\left(X,E\right)=\mathrm{Hom}_{\cur{N}\mathrm{IC}\left(X\right)}\left(\cur{O}_X,E\right)\end{equation}   and the adjunction becomes the identification
\begin{equation}
\mathrm{Hom}_{\cur{N}\mathrm{IC}\left(X\right)}\left(V\otimes_k\cur{O}_X,E\right)=\mathrm{Hom}_{\mathrm{Vec}_k}\left(V,\mathrm{Hom}_{\cur{N}\mathrm{IC}\left(X\right)}\left(\cur{O}_X,E\right)\right).
\end{equation}  
\end{example}

Since $f_*^\mathrm{dR}$ takes objects in $\cur{N}_f\mathrm{IC}\left(X\right)$ to objects in $\mathrm{IC}\left(S\right)$, it commutes with base change and there is an isomorphism of functors 
\begin{equation}H^0_\mathrm{dR}\left(X_s,-\right) \circ i_s^*\cong s^*\circ f_*^\mathrm{dR}:\cur{N}_f\mathrm{IC}\left(X\right)\rightarrow \mathrm{Vec}_k\end{equation}   (see for example \cite{Har75}, Chapter III, Theorem 5.2).

\begin{proposition} \label{trivfibiso} Suppose that $i_s^*E$ is trivial. Then the counit $\eps_{E}:f^*f_*^\mathrm{dR}E\rightarrow E$ is an isomorphism.
\end{proposition}

\begin{proof} Pulling back $\eps_E$ by $i_s^*$, and using base change, we get a morphism
\begin{equation} \cur{O}_{X_s}\otimes_k H^0_\mathrm{dR}\left(X_s,i_s^*E\right)\rightarrow i_s^*E
\end{equation}  
which by the explicit description of \ref{explad} is seen to be an isomorphism (as $i_s^*E$ is trivial). Hence by rigidity, $\eps_E$ must be an isomorphism. 
\end{proof}

\begin{proposition} \label{maxsubfib} Let $E\in\cur{N}_f\mathrm{IC}(X)$, and let $F_0\subset i_s^*(E)$ denote the largest trivial sub-object. Then there exists $E_0\subset E$ such that $F_0=i_s^*(E_0)$.
\end{proposition}

\begin{proof} Let $F=i_s^*(E)$. Since $H^0_\mathrm{dR}(X_s,F)=\mathrm{Hom}_{\mathrm{IC}(X_s)}(\cur{O}_{X_s},F)$, it follows that $F_0\cong\cur{O}_{X_s} \otimes_K H^0_\mathrm{dR}\left(X_s,F\right)$. Set $E_0=f^*f_*^{\mathrm{dR}}(E)$, then by the base change results proved above we know that $i_s^*(E_0)\cong F_0$, and that the natural map $E_0\rightarrow E$ restricts to the inclusion $F_0\rightarrow F$ on the fibre $X_s$. 
\end{proof}

\begin{corollary} \label{surject} The map $\pi_1^\mathrm{dR}\left(X_s,p_s\right)\rightarrow \pi_1^\mathrm{dR}\left(X/S,p\right)_s$ is a surjection.
\end{corollary}

We now turn to the proof of injectivity of the comparison map, borrowing heavily from ideas used in Section 2.1 of \cite{Had10}. We define objects $U_n$ of $\nic{X_s}$, the category of unipotent integrable connections on $X_s$ inductively as follows. $U_1$ will just be $\cur{O}_{X_s}$, and $U_{n+1}$ will be the extension of $U_{n}$ by $\cur{O}_{X_s}\otimes_k H^1_\mathrm{dR}\left(X_s,U_n^\vee\right)^\vee$ corresponding to the identity under the isomorphisms
\begin{align} \mathrm{Ext}_{\mathrm{IC}\left(X_s\right)}\left(U_n,\cur{O}_{X_s}\otimes_k H^1_\mathrm{dR}\left(X_s,U_n^\vee\right)^\vee\right) &
\cong H^1_\mathrm{dR}\left(X_s,U_n^\vee \otimes_k H^1_\mathrm{dR}\left(X_s,U_n^\vee\right)^\vee\right)  \\
&\cong H^1_\mathrm{dR}\left(X_s,U_n^\vee\right)\otimes_k H^1_\mathrm{dR}\left(X_s,U_n^\vee\right)^\vee  \\
&\cong \mathrm{End}_k\left(H^1_\mathrm{dR}\left(X_s,U_n^\vee\right)\right).
\end{align}
If we look at the long exact sequence in de Rham cohomology associated to the short exact sequence $0 \rightarrow U_n^\vee \rightarrow U_{n+1}^\vee \rightarrow H^1_\mathrm{dR}\left(X_s,U_n^\vee\right)\otimes_k\cur{O}_{X_s}\rightarrow 0$ we get 
\begin{align} 0\rightarrow H^0_\mathrm{dR}\left(X_s,U_n^\vee\right)&\rightarrow H^0_\mathrm{dR}\left(X_s,U_{n+1}^\vee\right)\rightarrow H^1_\mathrm{dR}\left(X_s,U_n^\vee\right)\\
  &
 \overset{\delta}{\rightarrow} H^1_\mathrm{dR}\left(X_s,U_n^\vee\right)\rightarrow H^1_\mathrm{dR}\left(X_s,U_{n+1}^\vee\right). \nonumber
\end{align}

\begin{lemma} The connecting homomorphism $\delta$ is the identity.
\end{lemma}

\begin{proof} By dualising, the extension 
\begin{equation} 0\rightarrow U_n^\vee \rightarrow U_{n+1}^\vee\rightarrow \cur{O}_{X_s}\otimes_k H^1_\mathrm{dR}\left(X_s,U_n^\vee\right)\rightarrow 0
\end{equation}  
corresponds to the identity under the isomorphism
\begin{equation}
\mathrm{Ext}_{\mathrm{IC}\left(X_s\right)}\left(\cur{O}_{X_s}\otimes_k H^1_\mathrm{dR}\left(X_s,U_n^\vee\right),U_n^\vee \right) \cong \mathrm{End}_k\left(H^1_\mathrm{dR}\left(X_s,U_n^\vee\right)\right)
\end{equation}  
Now the lemma follows from the fact that for an extension $0\rightarrow E\rightarrow F\rightarrow \cur{O}_{X_s}\otimes_k V\rightarrow 0$ of a trivial bundle by $E$, the class of the extension under the isomorphism 
\begin{equation} \mathrm{Ext}_{\mathrm{IC}\left(X_s\right)} \left(\cur{O}_{X_s}\otimes_k V,E\right)\cong V^\vee\otimes H^1_\mathrm{dR}\left(X_s,E\right)\cong \mathrm{Hom}_k\left(V,H^1_\mathrm{dR}\left(X_s,E\right)\right)
\end{equation}   is just the connecting homomorphism for the long exact sequence
\begin{equation} 0\rightarrow H^0_\mathrm{dR}\left(X_s,E\right)\rightarrow H^0_\mathrm{dR}\left(X_s,F\right)\rightarrow V\rightarrow H^1_\mathrm{dR}\left(X_s,E\right).
\end{equation}  
\end{proof}

In particular, any extension of $U_n$ by a trivial bundle $V\otimes_k \cur{O}_{X_s}$ is split after pulling back to $U_{n+1}$, and $H^0_\mathrm{dR}\left(X_s,U^\vee_{n+1}\right)\cong H^0_\mathrm{dR}\left(X_s,U^\vee_n\right)$. It then follows by induction that $H^0_\mathrm{dR}\left(X_s,U^\vee_n\right)\cong H^0_\mathrm{dR}\left(X_s,\cur{O}_{X_s}\right)\cong k$ for all $n$. 

\begin{definition} We define the unipotent class of an object $E\in\nic{X_s}$ inductively as follows. If $E$ is trivial, then we say $E$ has unipotent class 1. If there exists an extension 
\begin{equation} 0\rightarrow V\otimes_k \cur{O}_{X_s} \rightarrow E\rightarrow E'\rightarrow 0
\end{equation}  
with $E'$ of unipotent class $\leq m-1$, then we say that $E$ has unipotent class $\leq m$.
\end{definition}

Now let $x=p(s)$, $u_1=1\in \left(U_1\right)_x\cong \cur{O}_{X_s,x}=k$, and choose a compatible system of elements $u_n\in \left(U_n\right)_x$ mapping to $u_1$.

\begin{proposition} Let $F\in\nic{X_s}$ be an object of unipotent class $\leq m$. Then for all $n\geq m$ and any $f\in F_x$ there exists a morphism $\alpha:U_n\rightarrow F$ such that $\alpha_x\left(u_n\right)=f$.
\end{proposition}

\begin{proof} We copy the proof of Proposition 2.1.6 of \cite{Had10} and use induction on $m$. The case $m=1$ is straightforward. For the inductive step, let $F$ be of unipotent class $m$, and choose an exact sequence
\begin{equation} \label{uniextproof} 0\rightarrow E\overset{\psi}{\rightarrow}F\overset{\phi}{\rightarrow}G\rightarrow 0
\end{equation}  
with $E$ trivial and $G$ of unipotent class $<m$. By induction there exists a morphism $\beta:U_{n-1}\rightarrow G$ such that $\phi_x\left(f\right)=\beta_x\left(u_{n-1}\right)$. Pulling back the extension (\ref{uniextproof}) first by the morphism $\beta$ and then by the natural surjection $U_n\rightarrow U_{n-1}$ gives an extension of $U_n$ by $E$, which must split, as observed above.
\begin{equation}\xymatrix{
0 \ar[r]& E\ar[r]\ar@{=}[d]& F'' \ar[r]\ar[d]& U_n \ar[r]\ar[d]\ar@/_1pc/[l]& 0 \\
0 \ar[r]& E \ar[r]\ar@{=}[d]& F' \ar[r]\ar[d]& U_{n-1} \ar[r]\ar[d]& 0 \\
0 \ar[r]& E \ar[r]^\psi& F \ar[r]^\phi& G \ar[r]& 0
}
\end{equation}  
Let $\gamma:U_n\rightarrow F$ denote the induced morphism, then $\phi_x\left(\gamma_x\left(u_n\right)-f\right)=0$. Hence there exists some $e\in E_x$ such that $\psi_x\left(e\right)=\gamma_x\left(u_n\right)-f$. Again by induction we can choose $\gamma':U_n\rightarrow E$ with $\gamma'_x\left(u_n\right)=e$. Finally let $\alpha=\gamma - \psi\circ\gamma'$, it is easily seen that $\alpha_x(u_n)=f$. 
\end{proof}

\begin{corollary} Every $E$ in $\nic{X_s}$ is a quotient of $U_m^{\oplus N}$ for some $m,N\in \N$.
\end{corollary}

\begin{proof} Suppose that $E$ is of unipotent class $\leq m$. Let $e_1,\ldots ,e_N$ be a basis for $E_x$. Then there is a morphism $\alpha:U_m^{\oplus N}\rightarrow E$ with every $e_i$ in the image of the induced map on fibres. Thus $\alpha_x$ is surjective, and hence so is $\alpha$. 
\end{proof}

We now try to inductively define relatively nilpotent integrable connections $W_n$ on $X$ which restrict to the $U_n$ on fibres. Define higher direct images in de Rham cohomology by $\mathbf{R}^i_\mathrm{dR}f_*(E)=\cur{H}^{i-d}(f_+E)$, and begin the induction with $W_1=\cur{O}_X$. As part of the induction we will assume that $\RdR{0}\left(W_n^\vee\right)\cong \RdR{0}\left(\cur{O}_X\right)=\cur{O}_S$, that $\RdR{1}(W^\vee)$  and $\RdR{1}(W)$ are both coherent, i.e. regular integrable connections, and that there exists a horizontal morphism $p^*W_n\rightarrow\cur{O}_S$ such that the composite map $\cur{O}_S\cong g_*W_n^\vee\cong p^*g^*g_*W_n^\vee\rightarrow p^*W_n^\vee\rightarrow \cur{O}_S$ is an isomorphism. We will define $W_{n+1}$ to be an extension of $W_n$ by the sheaf $f^*\RdR{1}\left(W_n^\vee\right)^\vee$, and thus consider the extension group
\begin{equation} \mathrm{Ext}_{\mathrm{IC}\left(X\right)}\left(W_n,f^*\RdR{1}\left(W_n^\vee\right)^\vee\right)\cong H^1_\mathrm{dR}\left(X, W_n^\vee\otimes_{\cur{O}_X} f^*\RdR{1}\left(W_n^\vee\right)^\vee\right).
\end{equation}  

The Leray spectral sequence, together with the induction hypothesis and the projection formula, gives us the 4-term exact sequence
\begin{align} 0 &\rightarrow H^1_\mathrm{dR}\left(S,\RdR{1}\left(W_n^\vee\right)^\vee\right)\rightarrow \mathrm{Ext}_{\mathrm{IC}\left(X\right)}\left(W_n,f^*\RdR{1}\left(W_n^\vee\right)^\vee\right) \rightarrow \\ &\rightarrow \mathrm{End}_{\mathrm{IC}\left(S\right)}\left(\RdR{1}\left(W_n^\vee\right)\right)\rightarrow H^2_\mathrm{dR}\left(S,\RdR{1}\left(W_n^\vee\right)^\vee\right) \nonumber
\end{align}
and we can extract the commutative diagram
\begin{equation}
\xymatrix{
\mathrm{Ext}_{\mathrm{IC}\left(X\right)}\left(W_n,f^*\RdR{1}\left(W_n^\vee\right)^\vee\right) \ar[r]\ar[d] & \mathrm{Ext}_{\mathrm{IC}\left(X_s\right)}\left(U_n,\cur{O}_{X_s}\otimes_k H^1_\mathrm{dR}\left(X_s,U_n^\vee\right)^\vee\right) \ar@{=}[d] \\
\mathrm{End}_{\mathrm{IC}\left(S\right)}\left(\RdR{1}\left(W_n^\vee\right)\right) \ar[r]\ar[d] & \mathrm{End}_k\left(H^1_\mathrm{dR}\left(X_s,U_n^\vee\right) \right)\\
 H^2_\mathrm{dR}\left(S,\RdR{1}\left(W_n^\vee\right)^\vee\right)
}
\end{equation}  
where the horizontal arrows are just restrictions to fibres. The identity morphism in $\mathrm{End}_k\left(H^1_\mathrm{dR}\left(X_s,U_n^\vee\right)\right) $ clearly lifts to $\mathrm{End}_{\mathrm{IC}\left(S\right)}\left(\RdR{1}\left(W_n^\vee\right)\right)$, and hence the obstruction to finding $W_{n+1}$ lifting $U_{n+1}$ is the image of the identity under the map
\begin{equation} \mathrm{End}_{\mathrm{IC}\left(S\right)}\left(\RdR{1}\left(W_n^\vee\right)\right)\rightarrow H^2_\mathrm{dR}\left(S,\RdR{1}\left(W_n^\vee\right)^\vee\right).
\end{equation}  
In particular, if the base $S$ is an affine curve, this obstruction has to vanish.

\begin{proposition} Suppose $S$ is an affine curve. Then every object of $\nic{X_s}$ is a quotient of $\iota_s^*E$ for some $E\in \cur{N}_f\mathrm{IC}\left(X\right)$. 
\end{proposition}

\begin{proof}
To finish the induction step, we must show that 
\begin{equation}\RdR{0}\left(W_{n+1}^\vee\right)\cong\RdR{0}\left(W_n^\vee\right)\end{equation} and that $\RdR{1}(W_{n+1}^\vee)$ and $\RdR{1}(W_{n+1})$ are coherent. For the first claim, if we look at the long exact sequence of relative de Rham cohomology 
\begin{equation} 0 \rightarrow \RdR{0}\left(W_n^\vee\right) \rightarrow \RdR{0}\left(W_{n+1}^\vee\right)\rightarrow\ldots
\end{equation}  
we simply note that the given map restricts to an isomorphism on the fibre over $s$, and is hence an isomorphism. For the second, we simply use the long exact sequence in cohomology and the inductive hypothesis for $\RdR{1}\left(W_n\right)$ and $\RdR{1}\left(W_n^\vee\right)$.

\end{proof}

\begin{corollary} Suppose $S$ is an affine curve. Then the natural `base change' map $\pi_1^\mathrm{dR}(X_s,p_s)\rightarrow \pi_1^\mathrm{dR}(X/S,p)_s$ is an isomorphism.
\end{corollary}

\begin{remark} It is possible to define a relative fundamental group when $k$ is not necessarily algebraically closed (but still of characteristic $0$) using identical methods. One can then show that the corresponding `base change' question can be deduced from what we have proved in the algebraically closed case. Since this argument is rather fiddly, and not necessary in the context of this paper, we have omitted it.
\end{remark}

\section{Path torsors, non-abelian crystals and period maps}\label{tcp}

If $\cur{T}$ is a Tannakian category over an arbitrary field $k$, and $\omega_i$ are fibre functors on $\cur{T}$, $i=1,2$, with values in the category of quasi-coherent sheaves on some $k$-scheme $S$, then the functor of isomorphisms $\omega_1\rightarrow \omega_2$ is representable by an affine $S$-scheme, which is a $(G(\cur{T},\omega_1),G(\cur{T},\omega_2))$-bitorsor. This allows us to define path torsors under the algebraic and de Rham fundamental groups. In this section, we show how to do this in the relative case.

\subsection{Torsors in Tannakian categories}
\label{reltanna}

Let $\cur{C}$ be a Tannakian category over a field $k$. A Tannakian $\cur{C}$-category is a Tannakian category $\cur{D}$ together with an exact, $k$-linear tensor functor $t:\cur{C}\rightarrow \cur{D}$. We say it is neutral over $\cur{C}$ if there exists an exact, faithful $k$-linear tensor functor $\omega:\cur{D}\rightarrow \cur{C}$ such that $\omega\circ t\cong \mathrm{id}$. Such functors will be called fibre functors. If such a functor $\omega$ is fixed, we say $\cur{D}$ is neutralised. Thanks to \S6.10 of \cite{Del89}, we have a homomorphism
\begin{equation} t^*:\pi\left(\cur{D}\right)\rightarrow t\left(\pi\left(\cur{C}\right)\right)
\end{equation}  
of affine group schemes over $\cur{D}$. Hence applying $\omega$ gives us a homomorphism
\begin{equation} \omega\left(t^*\right):\omega\left(\pi\left(\cur{D}\right)\right)\rightarrow \pi\left(\cur{C}\right)
\end{equation}  
of affine group schemes over $\cur{C}$. We define $G\left(\cur{D},\omega\right):=\ker \omega\left(t^*\right)$.

For an affine group scheme $G$ over $\cur{C}$, let $\cur{O}_G$ be its Hopf algebra, a representation of $G$ is then defined to be an $\cur{O}_G$-comodule. That is an object $V\in \cur{C}$ together with a map $\delta:V\rightarrow \cur{O}_G \otimes V$ satisfying the usual axioms.

\begin{definition} \label{torsorsformal} A torsor under $G$ is a non-empty affine scheme $\mathrm{Sp}(\cur{O}_P)$ over $\cur{C}$, together with a $\cur{O}_G$-comodule structure on $\cur{O}_P$, such that the induced map $\cur{O}_P\otimes\cur{O}_P\rightarrow \cur{O}_P\otimes\cur{O}_G$ is an isomorphism.
\end{definition}

\begin{example} \label{actiontors} Suppose that $\cur{C}=\mathrm{Rep}_k\left(H\right)$, for some affine group scheme $H$ over $k$. Then an affine group scheme $G$ over $\cur{C}$ `is' just an affine group scheme $G_0$ over $k$ together with an action of $H$. A representation of $G$ `is' then just an $H$-equivariant representation of $G_0$, or in other words, a representation of the semi-direct product $G_0\rtimes H$.
\end{example}

Representations have another interpretation. Suppose that $V$ is an $\cur{O}_G$-comodule, and let $R$ be a $\cur{C}$-algebra. A point $g\in G\left(R\right)$ is then a morphism $\cur{O}_G\rightarrow R$ of $\cur{C}$-algebras, and hence for any such $g$ we get a morphism
\begin{equation} V\rightarrow V\otimes R
\end{equation}  
which extends linearly to a morphism
\begin{equation} V\otimes R\rightarrow V\otimes R.
\end{equation}  
This is an isomorphism, with inverse given by the map induced by $g^{-1}$. Hence we get an $R$-linear action of $G\left(R\right)$ on $V\otimes R$, for all $\cur{C}$-algebras $R$. The same proof as in the absolute case (Proposition 2.2 of \cite{MD81}) shows that a representation of $G$ (defined in terms of comodules) is equivalent to an $R$-linear action of $G\left(R\right)$ on $V\otimes R$, for all $R$.

For $G$ an affine group scheme over $\cur{C}$, let $\mathrm{Rep}_{\cur{C}}\left(G\right)$ denote its category of representations, this is a Tannakian category over $k$. There are canonical functors
\begin{equation}\xymatrix{\cur{C} \ar@<0.6ex>[r]^-{t} & \mathrm{Rep}_{\cur{C}}\left(G\right) \ar@<0.6ex>[l]^-{\omega}}
\end{equation}   given by `trivial representation' and `forget the representation'. This makes $\mathrm{Rep}_{\cur{C}}\left(G\right)$ neutral over $\cur{C}$. There is a natural homomorphism $G\rightarrow \omega(\pi(\mathrm{Rep}_\cur{C}(G)))$ which comes from the fact that by definition, $G$ acts on $\omega(V)$ for all $V\in\mathrm{Rep}_\cur{C}(G)$. Since this action is trivial on everything of the form $t(W)$, $W\in\cur{C}$, again by definition, this homomorphism factors to give a homomorphism
\begin{equation} G\rightarrow G(\mathrm{Rep}_{\cur{C}}(G),\omega).
\end{equation}
Conversely, if $\cur{D}$ is neutral over $\cur{C}$, with fibre functor $\omega$, then the action of $\omega(\pi(\cur{D}))$ on $\omega(V)$, for all $V\in\cur{D}$, induces an action of $G(\cur{D},\omega)$ on $\omega(V)$, and hence a functor
\begin{equation} \cur{D}\rightarrow \mathrm{Rep}_\cur{C}(G(\cur{D},\omega)).
\end{equation}

\begin{proposition} \label{relnt} In the above situation, the homomorphism 
\begin{equation}
G\rightarrow G(\mathrm{Rep}_{\cur{C}}(G),\omega)
\end{equation}
is an isomorphism, and the functor 
\begin{equation}\cur{D}\rightarrow \mathrm{Rep}_\cur{C}(G(\cur{D},\omega))
\end{equation} is an equivalence of categories. \end{proposition}

\begin{proof} If $\cur{C}$ is neutral, say $\cur{C}\cong \mathrm{Rep}_k(H)$, then thanks to Example \ref{actiontors} this is straightforward and amounts to little more than saying that the category of representations of a semi-direct product $G_0\rtimes H$ is equivalent to the category of $H$-equivariant $G_0$-representations. If $\cur{C}$ is not neutral, then we choose a fibre functor with values in some $k$-scheme $S$, apply Th\'{e}or\`{e}me 1.12 of \cite{Del90} and replace the affine group scheme $H$ by a certain groupoid acting on a $S$ (for more details see Section \ref{cohomperiod}). The argument is then formally identical.
\end{proof}

\begin{remark} Our definition of the fundamental group $\pi_1^\mathrm{dR}\left(X/S,p\right)$ is then just $G\left(\cur{N}_f\mathrm{IC}\left(X\right),
p^*\right)$, as an affine group scheme over $\mathrm{IC}\left(S\right)$. 
\end{remark}

In order to define torsors of isomorphisms in the relative setting, we must first recall Deligne's construction in the absolute case, which uses the notion of a coend. So suppose that we have categories $\cur{X}$ and $\cur{S}$, and a functor $F:\cur{X}\times \cur{X}^\mathrm{op}\rightarrow \cur{S}$. The coend of $F$ is the universal pair $\left(\zeta,s\right)$ where $s$ is an object of $\cur{S}$ and $\zeta:F\rightarrow s$ is a bi-natural transformation. Here $s$ is the constant unctor at $s\in \mathrm{Ob}\left(\cur{S}\right)$, and by bi-natural we mean that it is natural in both variables. If such an object exists, we will denote it by
\begin{equation} \int^\cur{X} F\left(x,x\right).
\end{equation}  
If $\cur{S}$ is cocomplete then the coend always exists and is given concretely by the formula (see Chapter IX, Section 6 of \cite{Mac71})
\begin{equation} \int^\cur{X} F\left(x,x\right)= \mathrm{colim}\left( \CMcoprod_{f:x\rightarrow y\in \mathrm{Mor}\left(\cur{X}\right)}F\left(x,y\right) \rightrightarrows \CMcoprod_{x\in\mathrm{Ob}\left(\cur{X}\right)}F\left(x,x\right) \right).
\end{equation}  

Suppose that $\cur{C}$ is a Tannakian category, and let $\omega_1,\omega_2:\cur{C}\rightarrow \mathrm{Qcoh}\left(S\right)$ be two fibre functors on $\cur{C}$. In \cite{Del90}, Deligne defines
\begin{equation} L_S\left(\omega_1,\omega_2\right)=\int^\cur{C} \omega_1\left(V\right) \otimes \omega_2\left(V\right)^\vee	
\end{equation}  
to be the coend of the bifunctor
\begin{equation} \omega_1\otimes \omega_2^\vee : \cur{C}\times \cur{C}^\mathrm{op}\rightarrow \mathrm{Qcoh}\left(S\right),
\end{equation}  
and in \S6 of \emph{loc. cit.}, uses the tensor structure of $\cur{C}$ to define a multiplication on $L_S\left(\omega_1,\omega_2\right)$ which makes it into a quasi-coherent $\cur{O}_S$-algebra. He then proves that $\mathrm{Spec}\left(L_S\left(\omega_1,\omega_2\right)\right)$ represents the functor of isomorphisms from $\omega_1$ to $\omega_2$. 

Now let $\cur{C}$ be a Tannakian category, let $\cur{D}$ be neutral over $\cur{C}$, and suppose that $\omega_1,\omega_2:\cur{D}\rightarrow \cur{C}$ are two fibre functors from $\cur{D}$ to $\cur{C}$. Define the coend
\begin{equation} L_\cur{C}\left(\omega_1,\omega_2\right):=\int^\cur{D}\omega_1\left(V\right)\otimes \omega_2\left(V\right)^\vee \in\mathrm{Ind}\left(\cur{C}\right).
\end{equation}   
If $\eta:\cur{C}\rightarrow \mathrm{Qcoh}\left(S\right)$ is a fibre functor, then $\eta$ commutes with colimits, and hence $\eta\left(L_\cur{C}\left(\omega_1,\omega_2\right)\right)=L_S\left(\eta\omega_1,\eta\omega_2\right)$: this is a quasi-coherent $\cur{O}_S$-algebra, functorial in $\eta$. Since algebraic structures in Tannakian categories, such as commutative algebras, Hopf algebras, and so on, can be constructed `functorially in fibre functors', (see for example \S5.11 of \cite{Del89}), it follows that there is a unique way of defining a $\cur{C}$-algebra structure on $L_\cur{C}\left(\omega_1,\omega_2\right)$ lifting the $\cur{O}_S$-algebra structure on each $\eta\left(L_\cur{C}\left(\omega_1,\omega_2\right)\right)$. Moreover, since $\eta\left(\mathrm{Sp}\left(L_\cur{C}\left(\omega_1,\omega_2\right)\right)\right)$ is  a $\left(\eta\omega_1\left(\pi\left(\cur{D}\right)\right),\eta\omega_2\left(\pi\left(\cur{D}\right)\right)\right)$-bitorsor, functorially in $\eta$, the affine scheme 
\begin{equation}
P_\cur{C}\left(\omega_1,\omega_2\right):=\mathrm{Sp}\left(L_\cur{C}\left(\omega_1,\omega_2\right)\right)\end{equation}   is a $\left(\omega_1\left(\pi\left(\cur{D}\right)\right),\omega_2\left(\pi\left(\cur{D}\right)\right)\right)$-bitorsor over $\cur{C}$.

What we actually want, however, is a $\left(G_{\cur{C}}\left(\cur{D},\omega_2\right),G_{\cur{C}}\left(\cur{D},\omega_2\right)\right)$-bitorsor. We get this as follows. Suppose that $V\in\cur{D}$, then by the definition of $L_\cur{C}\left(\omega_1,\omega_2\right)$ we get a morphism
\begin{equation} \omega_1\left(V\right) \otimes \omega_2\left(V\right)^\vee\rightarrow L_\cur{C}\left(\omega_1,\omega_2\right)
\end{equation}  
which corresponds to a morphism
\begin{equation} \omega_1\left(V\right)\rightarrow \omega_2\left(V\right)\otimes L_\cur{C}\left(\omega_1,\omega_2\right).
\end{equation}  
Thus a morphism $L_\cur{C}\left(\omega_1,\omega_2\right)\rightarrow R$ for some $\cur{C}$-algebra $R$ induces an $R$-linear morphism
\begin{equation} \omega_1\left(V\right)\otimes R\rightarrow \omega_2\left(V\right)\otimes R
\end{equation}  
which is in fact an isomorphism, since it is so after applying any fibre functor.

\begin{definition} Define $P_\mathrm{triv}\left(\omega_1,\omega_2\right)$ to be the sub-functor of $P_\cur{C}\left(\omega_1,\omega_2\right)$ which takes $R$ to the set of all morphisms $L_\cur{C}\left(\omega_1,\omega_2\right)\rightarrow R$ such that for every $V$ in the essential image of $t:\cur{C}\rightarrow \cur{D}$, the induced automorphism of $R\otimes \omega_1(V)=R\otimes\omega_2(V)$ is the identity.
\end{definition}

\begin{proposition} The functor $P_\mathrm{triv}\left(\omega_1,\omega_2\right)$ is representable by an affine scheme over $\cur{C}$, and is a $\left(G_{\cur{C}}\left(\mathcal{D},\omega_1\right),G_{\cur{C}}\left(\mathcal{D},\omega_2\right)\right)$-bitorsor in the category of affine schemes over $\mathcal{C}$.
\end{proposition}

\begin{proof}

First note that if $V\in\mathrm{Ob}\left(\mathcal{D}\right)$, then $\omega_i\left(\pi\left(\mathcal{D}\right)\right)$ acts on $\omega_i\left(V\right)$, and $G\left(\cur{D},\omega_i\right)$ is the largest subgroup of $\omega_i\left(\pi\left(\mathcal{D}\right)\right)$ whose action on $\omega_i\left(V\right)$ is trivial for all $V$ in the essential image of $t$.

Now, if $p\in P_\mathrm{triv}\left(\omega_1,\omega_2\right)\left(R\right)$ and $g\in G_{\cur{C}}\left(\cur{D},\omega_1\right)\left(R\right)$ then $gp\in P_\cur{C}\left(\omega_1,\omega_2\right)\left(R\right)$ acts trivially on everything of the form $t\left(W\right)$, and hence lies in $P_\mathrm{triv}\left(\omega_1,\omega_2\right)\left(R\right)$. Hence $G\left(\cur{D},\omega_1\right)$ acts on $P_\mathrm{triv}\left(\omega_1,\omega_2\right)$. For $p,p'\in P_\mathrm{triv}\left(\omega_1,\omega_2\right)\left(R\right)$, $p^{-1}p'$ is an automorphism of $\omega_1(V)\otimes R$ which is trivial for all $V$ in the essential image of $t$. Hence it must be an element of $G\left(\cur{D},\omega_1\right)\left(R\right)\subset \omega\left(\pi_1\left(\cur{D}\right)\right)\left(R\right)$. The same arguments work for $G_{\cur{C}}\left(\cur{D},\omega_2\right)$.

Thus $P_\mathrm{triv}(\omega_1,\omega_2)$ is a bi-pseudo-torsor, and to complete the proof, we must show that $P_\mathrm{triv}\left(\omega_1,\omega_2\right)$ is represented by a non-empty affine scheme over $\cur{C}$. By similar arguments to before, one can see that the fundamental group $\pi\left(\cur{C}\right)$ of $\cur{C}$ is the formal $\mathrm{Spec}$ of the Hopf $\cur{C}$-algebra
\begin{equation} L_\cur{C}\left(\mathrm{id},\mathrm{id}\right)=\mathrm{colim}\left( \CMcoprod_{f:V\rightarrow W\in \mathrm{Mor}\left(\cur{C}\right)}V\otimes W^\vee \rightrightarrows \CMcoprod_{V\in\mathrm{Ob}\left(\cur{C}\right)}V\otimes V^\vee \right)
\end{equation}  
and hence one can construct a morphism of affine $\cur{C}$-schemes
\begin{equation} P_\cur{C}\left(\omega_1,\omega_2\right)\rightarrow \pi\left(\cur{C}\right)
\end{equation}  
which is the formal $\mathrm{Spec}$ of the obvious morphism 
$ L_\cur{C}\left(\mathrm{id},\mathrm{id}\right)\rightarrow  L_\cur{C}\left(\omega_1,\omega_2\right). $
Then $P_\mathrm{triv}\left(\omega_1,\omega_2\right)$ is the fibre of $P_\cur{C}\left(\omega_1,\omega_2\right)\rightarrow \pi\left(\cur{C}\right)$ over the identity section $\mathrm{Sp}\left(1\right)\rightarrow \pi\left(\cur{C}\right)$. Hence it is the formal $\mathrm{Spec}$ of the algebra $L_\mathrm{triv}\left(\omega_1,\omega_2\right)$ defined by the push-out diagram
\begin{equation}
\xymatrix{
L_\cur{C}\left(\mathrm{id},\mathrm{id}\right)\ar[r]\ar[d] & 1\ar[d] \\
L_\cur{C}\left(\omega_1,\omega_2\right) \ar[r] & L_\mathrm{triv}\left(\omega_1,\omega_2\right)
}
\end{equation}  
and is thus representable by an affine $\cur{C}$-scheme.

To prove that $P_\mathrm{triv}\left(\omega_1,\omega_2\right)\neq\emptyset$, it suffices to show that $\eta\left(P_\mathrm{triv}\left(\omega_1,\omega_2\right)\right)\neq\emptyset$ for any fibre functor $\eta:\cur{C}\rightarrow \mathrm{Qcoh}\left(S\right)$. For any $f:T\rightarrow S$, $\eta\left(P_\mathrm{triv}\left(\omega_1,\omega_2\right)\right)\left(T\right)$ is the subset of $\mathrm{Isom}^\otimes\left(f^*\circ\eta\omega_1,f^*\circ\eta\omega_2\right)$ which maps to the identity under the natural map
\begin{equation} r: \mathrm{Isom}^\otimes\left(f^*\circ\eta\omega_1,f^*\circ\eta\omega_2\right)\rightarrow \mathrm{Isom}^\otimes(f^*\circ\eta\omega_1t,f^*\circ\eta\omega_2t)=\mathrm{Aut}^\otimes \left(f^*\circ \eta\right).
\end{equation}  
There is certainly some $S$-scheme $f:T\rightarrow S$ such that the LHS is non-empty. Pick such a $T$, and pick some $p\in \mathrm{Isom}^\otimes\left(f^*\circ\eta\omega_1,f^*\circ\eta\omega_2\right)$. Since the morphism $\omega_1\left(\pi\left(\cur{D}\right)\right)\rightarrow \pi\left(\cur{C}\right)$ admits a section, the induced homomorphism
\begin{equation} \mathrm{Aut}^\otimes \left(f^*\circ \eta\omega_1\right) \rightarrow \mathrm{Aut}^\otimes \left(f^*\circ \eta\right)
\end{equation}  
is surjective, and hence there exists some $g\in \mathrm{Aut}^\otimes \left(f^*\circ \eta\omega_1\right)$ mapping to $r\left(p\right)\in \mathrm{Aut}^\otimes \left(f^*\circ \eta\right)$. Then $p':=g^{-1}p$ is an element of the set $\mathrm{Isom}^\otimes\left(f^*\circ\eta\omega_1,f^*\circ\eta\omega_2\right)$ and $r\left(p'\right)=\mathrm{id}$, thus $\eta\left(P_\mathrm{triv}\left(\omega_1,\omega_2\right)\right)\left(T\right)\neq\emptyset$.
\end{proof}

\begin{remark} We can rephrase this as follows. Consider the functors of $\cur{C}$ algebras
\begin{align} \underline{\mathrm{Isom}}^\otimes\left(\omega_1,\omega_2\right):\cur{C}\mathrm{-alg}&\rightarrow \left(\mathrm{Set}\right) \nonumber\\ R &\mapsto \mathrm{Isom}^\otimes\left(\omega_1\left(-\right)\otimes R,\omega_2\left(-\right)\otimes R\right); \\ 
\underline{\mathrm{Aut}}^\otimes\left(\mathrm{id}\right):\cur{C}\mathrm{-alg}&\rightarrow \left(\mathrm{Set}\right) \nonumber\\ R &\mapsto \mathrm{Aut}^\otimes\left(\left(-\right)\otimes R\right);  
\end{align}
as well as the sub-functor  $\underline{\mathrm{Isom}}_\cur{C}^\otimes\left(\omega_1,\omega_2\right)$,
the `functor of $\cur{C}$-isomorphisms $\omega_1\rightarrow \omega_2$', defined to be the fibre over the identity of the natural morphism
\begin{equation} \underline{\mathrm{Isom}}^\otimes\left(\omega_1,\omega_2\right) \rightarrow \underline{\mathrm{Aut}}^\otimes\left(\mathrm{id}\right).
\end{equation}  
Then the functor $\underline{\mathrm{Isom}}_\cur{C}^\otimes\left(\omega_1,\omega_2\right)$ is representable by the affine scheme $P_\mathrm{triv}\left(\omega_1,\omega_2\right)$ over $\cur{C}$, which is a $\left(G_{\cur{C}}\left(\mathcal{D},\omega_1\right),G_{\cur{C}}\left(\mathcal{D},\omega_2\right)\right)$ bitorsor. 
\end{remark}

\subsection{Path torsors under relative fundamental groups}

Let $k$ be an algebraically closed field of characteristic zero, $S$ a connected, affine curve over $k$ and $f:X\rightarrow S$ a `good' morphism. Let $p,x$ be sections of $f$. 
We can apply the above methods to obtain an affine scheme over $\mathrm{IC}\left(S\right)$, the torsor of paths from $x$ to $p$, which can be considered as an affine scheme $P\left(x\right)=\pi_1^\mathrm{dR}\left(X/S,x,p\right)$ over $S$, together with an integrable connection on $\cur{O}_{P\left(x\right)}$ (as a quasi-coherent $\cur{O}_S$ algebra). This is naturally a left torsor under $\pi_1^\mathrm{dR}\left(X/S,x\right)$ and a right torsor under $\pi_1^\mathrm{dR}\left(X/S,p\right)=:G$. Moreover, the action map $P\left(x\right)\times G\rightarrow P\left(x\right)$ is compatible with the connections, in the sense that the associated comodule structure
\begin{equation} \cur{O}_{P\left(x\right)}\rightarrow \cur{O}_{P\left(x\right)}\otimes_{\cur{O}_S} \cur{O}_G
\end{equation}  
is horizontal, the RHS being given the tensor product connection. If $G_n$ is the quotient of $G$ by the $n$th term in its lower central series, we will denote the push-out torsor $P\left(x\right)\times^G G_n$ by $P\left(x\right)_n$. As before, the action map $P\left(x\right)_n\times G_n\rightarrow P\left(x\right)_n$ is compatible with the connections.

\begin{definition} A $\nabla$-torsor under $G_n$ is a $G_n$-torsor $P$ over $S$ in the usual sense, together with a regular integrable connection on $\cur{O}_P$, such that the action map 
\begin{equation} \cur{O}_P\rightarrow \cur{O}_P\otimes\cur{O}_{G_n}
\end{equation}
is horizontal. The set of isomorphism classes of $\nabla$-torsors is denoted $H^1_\nabla(S,G_n)$.
\end{definition}

Thus we have `period maps'
\begin{equation}X(S)\rightarrow H^1_\nabla(S,G_n)
\end{equation}
which takes $x\in X(S)$ to the path torsor $P(x)_n$.

\begin{remark}
\begin{enumerate}
\item This is not a good period map to study. For instance, if $k=\C$, then the relative fundamental group is not just an affine group scheme with connection. There are reasons to expect that one can put a `non-abelian' variation of Hodge structure on this fundamental group. Similar considerations will apply to the path torsors, and the period maps should take these variations of Hodge structures into account. 
\item We can use the pro-nilpotent Lie algebra of $\pi_1^\mathrm{dR}(X/S,p)$ and the Campell-Hausdorff law to view $\pi_1^\mathrm{dR}(X/S,p)$ as a non-abelian sheaf of groups on the infinitesimal site of $S/k$. We can use this interpretation to give an alternative definition of the cohomology set $H^1_\nabla(S,G_n)$.
\item A natural question to ask is whether or not, as in the situation studied by Kim, the targets for the period maps have the structure of algebraic varieties. Since we are more interested in the positive characteristic case, we will not pursue this question here.
\end{enumerate}
\end{remark}

\section{Crystalline fundamental groups of smooth families in char $p$}

Our goal in this chapter is to define the fundamental group of a smooth family $f:X\rightarrow S$ of varieties over a finite field. Many of our arguments are essentially the same as those we gave in Chapter \ref{zero}

We will assume that the reader is familiar with the theory of rigid cohomology and overconvergent ($F$-)isocrystals, a good reference is \cite{Ber96b}. Assume that $k$ is a finite field, of order $q=p^a$ and characteristic $p>0$. Frobenius will always refer to linear Frobenius. If $U/K$ is a variety, the category  of overconvergent ($F$-)isocrystals on $U/K$ is denoted $(F\text{-})\mathrm{Isoc}^\dagger(U/K)$. These are Tannakian categories over $K$.

We define $\cur{N}\mathrm{Isoc}^\dagger(U/K)$ to be the full subcategory of $\mathrm{Isoc}^\dagger(U/K)$ on objects admitting a filtration whose graded pieces are constant. Chiarellotto and Le Stum in \cite{CLS99a} define the rigid fundamental group $\pi_1^\text{rig}(U,x)$ of $U$ at a $k$-rational point $x$ to be the Tannaka dual of $\cur{N}\text{Isoc}^\dagger(U/K)$ with respect to the fibre functor $x^*$. This is a pro-unipotent group scheme over $K$.

Now suppose that $g:X\rightarrow S$ is a `good', proper morphism over $k$, and let $p:S\rightarrow X$ be a section.

\begin{definition} We say that $E\in F\text{-}\text{Isoc}^\dagger(X/K)$ is relatively unipotent if there is a filtration of $E$, whose graded pieces are all in the essential image of $g^*:F\text{-}\text{Isoc}^\dagger(S/K)\rightarrow F\text{-}\text{Isoc}^\dagger(X/K)$. The full subcategory of relatively unipotent overconvergent $F$-isocrystals is denoted $\cur{N}_gF\text{-}\text{Isoc}^\dagger(X/K)$. 
\end{definition} 

The pair of functors
\begin{equation}\xymatrix{ \cur{N}_gF\text{-}\text{Isoc}^\dagger(X/K) \ar@<0.6ex>[r]^-{p^*} & \ar@<0.6ex>[l]^-{g^*} F\text{-}\text{Isoc}^\dagger(S/K)}
\end{equation}  
makes $\cur{N}_gF\text{-}\text{Isoc}^\dagger(X/K)$ neutral over $F\text{-}\text{Isoc}^\dagger(S/K)$ in the sense of \S2.1. Hence we get an affine group scheme $G(\cur{N}_gF\text{-}\mathrm{Isoc}^\dagger(X/K),p^*)$ in $F\text{-}\text{Isoc}^\dagger(S/K)$.

\begin{definition} We define the relative fundamental group to be the affine group scheme $G(\cur{N}_gF\text{-}\mathrm{Isoc}^\dagger(X/K),p^*)$ in $F\text{-}\text{Isoc}^\dagger(S/K)$.
\end{definition}

For $s\in S$ a closed point, let $i_s:X_s\rightarrow X$ denote the inclusion of the fibre over $s$ and  let $g_s:X_s\rightarrow \spec{k(s)}$ denote the structure morphism. Let $K(s)$ denote the unique unramified extension of $K$ with residue field $k(s)$. Let $\cur{V}(s)$ denote the ring of integers of $K(s)$. In keeping with notation of previous chapters, let $\pi_1^\text{rig}(X/S,p)_s$ denote the affine group scheme $s^*(\pi_1^\text{rig}(X/S,p))$ over $K(s)$. The pull-back functor 
\begin{equation} i_s^*:\cur{N}_gF\text{-}\text{Isoc}^\dagger(X/K)\rightarrow \cur{N}\text{Isoc}^\dagger(X_s/K(s))\end{equation}   induces a homomorphism
\begin{equation} \phi:\pi_1^\text{rig}(X_s,p_s)\rightarrow \pi_1^\text{rig}(X/S,p)_s
\end{equation}  
of affine group schemes over $K$. We would like to show again that when $S$ is an affine curve, this is an isomorphism. The question is whether or not the sequence of affine group schemes corresponding to the sequence of neutral Tannakian categories
\begin{equation} \cur{N}\mathrm{Isoc}(X_s/K(s)) \leftarrow \cur{N}_gF\text{-}\text{Isoc}^\dagger(X/K)\otimes_K K(s) \leftarrow F\text{-}\text{Isoc}^\dagger(S/K)\otimes_K K(s)
\end{equation}  
is exact. Thus, as before, this boils down to the following three questions.

\begin{enumerate}\item If $E\in\cur{N}_gF\text{-}\text{Isoc}^\dagger(X/K)\otimes_K K(s)$ is such that $i_s^*E$ is constant, is $E$ of the form $g^*F$ for some $F\in F\text{-}\text{Isoc}^\dagger(S/K)\otimes_K K(s)$?
\item If $E\in \cur{N}_gF\text{-}\mathrm{Isoc}^\dagger(X/K)\otimes_K K(s)$, and $F_0\subset i_s^*E$ denotes the largest constant sub-object, then does there exist $E_0\subset E$ such that $F_0=i_s^*E_0$?
\item Given $E\in \mathrm{Isoc}^\dagger(X_s/K(s))$, does there exist $F\in\cur{N}_gF\text{-}\text{Isoc}^\dagger(X/K)\otimes_K K(s)$ such that $E$ is a quotient of $i_s^*F$?
\end{enumerate}

\begin{remark} Actually, in order to apply these criteria, we need to know that the kernel of the homomorphism of group schemes corresponding to 
\begin{equation} \cur{N}_gF\text{-}\text{Isoc}^\dagger(X/K)\otimes_K K(s) \leftarrow F\text{-}\text{Isoc}^\dagger(S/K)\otimes_K K(s)
\end{equation}  
is pro-unipotent, or using Lemma 1.3, Part I of \cite{Wil97}, that every object $E$ of the category $\cur{N}_gF\text{-}\text{Isoc}^\dagger(X/K)\otimes_K K(s)$ has a non-zero subobject of the form $f^*F$ for some $F\in F\text{-}\text{Isoc}^\dagger(S/K)\otimes_K K(s)$. Let $E_0$ denote the largest relatively constant sub-object of $E$, considered in the category $\cur{N}_gF\text{-}\text{Isoc}^\dagger(X/K)$. Then functoriality of $E_0$ implies that a $K(s)$ module structure $K(s)\rightarrow\mathrm{End}(E)$ will induce one on $E_0$. Hence we must show that an $K(s)$-module structure on $f^*F$ induces one on $F$. But now just use the section $p$ to get a homomorphism of rings $\mathrm{End}(f^*F)\rightarrow \mathrm{End}(F)$. 
\end{remark}

As in the case of characteristic $0$, we will only show that the base change map is an isomorphism when the base is an affine curve, and under some mild technical hypotheses on $X$. We will then use a gluing argument to construct $\pi_1^\rig(X/S,p)$ for (not necessarily affine) curves.

\subsection{Base change for affine curves}

Hypotheses and notations will be as in the previous section, except that we now assume that $S$ is a smooth affine curve. We will make the following additional technical hypothesis.

\begin{hypothesis} \label{hypothesis} There exists a smooth and proper formal $\cur{V}$-scheme $\cur{P}$, an immersion $X\rightarrow P$ of $X$ into its special fibre, such that the closure $X'$ of $X$ in $P$ is smooth, and there exists a divisor $T$ of $P$ with $X=X'\setminus T$. 
\end{hypothesis}

\begin{remark} \begin{enumerate}\item We should eventually be able to remove this technical hypothesis, using methods of `recollement', but we do not worry about this for now.
\item One non-trivial example of such a $g$ is given by a model for a smooth, proper, geometrically connected curve $C$ over a function field $K$ over a finite field. In this situation $S'$ is the unique smooth, proper model for $K$, $X'$ is a regular, flat, proper $S'$-scheme, whose generic fibre is $C$, $S\subset S'$ is an affine open subset of $S'$ over which $g$ is smooth, and $X$ is the pre-image of $S$. Since $X'$ is a regular, proper surface over a finite field, it is smooth, hence projective, and the above hypotheses really are satisfied.
\item Since $S$ is a smooth curve, these technical hypotheses are automatically satisfied for $S$. 
\end{enumerate}
\end{remark}

In this section we will prove the following two theorems.

\begin{theorem} \label{one} \begin{enumerate} \item Let $E\in \cur{N}_{g}F\text{-}\mathrm{Isoc}^\dagger(X/K)\otimes_K K(s) $ and suppose that $i_s^*E$ is a constant isocrystal. Then there exists $E'\in F\text{-}\mathrm{Isoc}^\dagger(S/K)\otimes_K K(s)$ such that $E\cong g^*E'$.
\item Let $E\in \cur{N}_{g}F\text{-}\mathrm{Isoc}^\dagger(X/K)\otimes_K K(s) $, and let $F_0\subset i_s^*E$ denote the largest constant subobject. Then there exists $E_0\subset E$ such that $F_0=i_s^*E_0$.
\end{enumerate}
\end{theorem}

\begin{theorem} \label{two} Let $E\in \cur{N}\mathrm{Isoc}^\dagger(X_s/K(s))$. Then there exists some object $E'\in \cur{N}_{g}F\text{-}\mathrm{Isoc}^\dagger(X/K)\otimes_K K(s)$ such that $E$ is a quotient of $i_s^*E'$.
\end{theorem}

\begin{remark} The reason we have used categories of overconvergent $F$-isocrystals rather than  overconvergent isocrystals without Frobenius is that the theory of `six operations' has only fully been developed for overconvergent $F$-isocrystals. If six operations were to be resolved for overconvergent isocrystals in general, then we would be able to deduce results for smooth fibrations over any perfect field of positive characteristic, not just over finite fields where we can linearise Frobenius.
\end{remark}

The method of proof will be entirely analogous to the proof in characteristic $0$, replacing the algebraic $\cur{D}$-modules used there by their arithmetic counterparts, the theory of which was developed by Berthelot and Caro. It would be far too much of a detour to describe this theory in any depth, so instead we will just recall the notations and results needed, referring the reader to the series of articles \cite{Ber02}, \cite{Ber96a}, \cite{Ber00} and \cite{Car15b}, \cite{Car09b}, \cite{Car04}, \cite{Car15a}, \cite{Car07}, \cite{Car06a} for details.

We let $F\text{-}D^b_\mathrm{surhol}(\cur{D}_{X/K})$ (resp. $F\text{-}D^b_\mathrm{surhol}(\cur{D}_{S/K})$) denote the category of overholonomic $F\text{-}\cur{D}$-modules on $X$ (resp. $S$) as defined in Section 3 of \cite{Car09b}. There is a functor 
\begin{equation} \mathrm{sp}_{X,+}:F\text{-}\mathrm{Isoc}^\dagger(X/K)\rightarrow F\text{-}D^b_\mathrm{surhol}(\cur{D}_{X/K})
\end{equation}
which is an equivalence onto the full subcategory $F\text{-}\mathrm{Isoc}^{\dagger\dagger}(X/K)$ of overcoherent $F$-isocrystals (Theorem 2.3.16 of \cite{CT12} and Th\'{e}or\`{e}me 2.3.1 of \cite{Car07}) and compatible with the natural tensor products on both sides (Proposition 4.8 of \cite{Car15b}). The same also holds for $S$. Let 
\begin{align} g_+:F\text{-}D^b_\mathrm{surhol}(\cur{D}_{X/K})&\rightarrow F\text{-}D^b_\mathrm{surhol}(\cur{D}_{S/K}) \\
g^+:F\text{-}D^b_\mathrm{surhol}(\cur{D}_{S/K})&\rightarrow F\text{-}D^b_\mathrm{surhol}(\cur{D}_{X/K})
\end{align} be the adjoint functors defined in Section 3 of \cite{Car09b}. By Th\'{e}or\`{e}me 4.2.12 of \cite{Car15a}, for any $E\in F\text{-}\mathrm{Isoc}^\dagger(X/K)$, and any $i\in\Z$, $\cur{H}^i(g_+\mathrm{sp}_{X,+}(E))\in F\text{-}\mathrm{Isoc}^{\dagger\dagger}(S/K)$ and hence we can define
\begin{equation} g_*:=sp_{S,+}^{-1}\cur{H}^{-d}(g_+\mathrm{sp}_{X,+}(-))(-d):F\text{-}\mathrm{Isoc}^{\dagger}(X/K)\rightarrow  F\text{-}\mathrm{Isoc}^{\dagger}(S/K)
\end{equation} 
where $d$ is the relative dimension of $X/S$, and $(-d)$ denotes the Tate twist. We can also define the higher direct images
\begin{equation} \mathbf{R}^ig_*:=sp_{S,+}^{-1}\cur{H}^{-d+i}(g_+\mathrm{sp}_{X,+}(-))(-d):F\text{-}\mathrm{Isoc}^{\dagger}(X/K)\rightarrow  F\text{-}\mathrm{Isoc}^{\dagger}(S/K).
\end{equation} 
Let $s^!:F\text{-}D^b_\mathrm{surhol}(\cur{D}_{S/K})\rightarrow F\text{-}D^b_\mathrm{surhol}(\cur{D}_{\spec{k(s)}/K(s)})$ denote the functor defined in Section 3 of \cite{Car09b}. 

\begin{remark} Although Caro's functor $s^!$ lands in $F\text{-}D^b_\mathrm{surhol}(\cur{D}_{\spec{k(s)}/K})$ rather than $F\text{-}D^b_\mathrm{surhol}(\cur{D}_{\spec{k(s)}/K(s)})$, it can be easily adapted to land in the latter category. The base change result that we use below holds in this slightly altered context.
\end{remark}

\begin{proposition}\label{newbasechange} Let $s\in S$ be a closed point. There is an isomorphism of functors $s^*\mathbf{R}^ig_*(-)\cong H^i_\mathrm{rig}(X_s,i_s^*(-)):F\text{-}\mathrm{Isoc}^\dagger(X/K)\rightarrow \mathrm{Vec}_{K(s)}$.
\end{proposition}

\begin{remark} We are deliberately ignoring Frobenius structure in the final target category of these two composite functors.
\end{remark}

\begin{proof} This follow from proper base change for arithmetic $\cur{D}$-modules (Th\'{e}or\`{e}me 4.4.2 of \cite{Car15a}), together with the identification $s^*=s^![1]$ for overcoherent $F$-isocrystals on $S$ (1.4.5 of \cite{Car15b}, recall $\dim S=1$) and the fact that, defining $\mathbf{R}^ig_{s*}$ entirely analogously to $g_*$, we have the identification $\mathbf{R}^ig_{s*}(-)=H^i_\mathrm{rig}(X_s,-)$ (since we are not worried about the Frobenius structure on $H^i_\mathrm{rig}(X_s,-)$, this follows from Lemme 7.3.4 of \cite{Car06a}).
\end{proof}

\begin{proposition} For $E\in F\text{-}\mathrm{Isoc}^\dagger(X/K)$, $g_+\mathrm{sp}_{X,+}(E)$ is concentrated in degrees $\geq -d$.
\end{proposition}

\begin{proof} We know that $g_+\mathrm{sp}_{X,+}(E)$ has overcoherent $F$-isocrystals for cohomology sheaves, and by the previous proposition, the fibre over $s$ of $\cur{H}^i(g_+\mathrm{sp}_{X,+}(E))$ is zero for $i\leq -d$. Hence $\cur{H}^i(g_+\mathrm{sp}_{X,+}(E))$ is zero for $i\leq -d$.
\end{proof}

\begin{proposition} $g_*$ is right adjoint to $g^*$.
\end{proposition}

\begin{proof} Since $g_+$ is right adjoint to $g^+$, this just follows from the previous proposition and the fact that $g^+\mathrm{sp}_{S,+}(-)[d](d)=\mathrm{sp}_{X,+}g^*(-)$.
\end{proof}

\begin{proof}[Proof of Theorem \ref{one}] Because $g_*$ and $g^*$ are functorial, they extend to give adjoint functors
\begin{align} \xymatrix{g^*:F\text{-}\mathrm{Isoc}^\dagger(S/K)\otimes_K K(s)\ar@<0.6ex>[r]& \ar@<0.6ex>[l] \cur{N}_{g}F\text{-}\mathrm{Isoc}^\dagger\left(X/K\right)\otimes_K K(s):{g}_*}
\end{align}  
such that (using the base change theorem as in the proof of Proposition \ref{newbasechange}) the counit $g^*g_*E\rightarrow E$ restricts to the counit of the adjunction
\begin{equation} \xymatrix{-
\otimes_{K(s)} \cur{O}_{X_s/K(s)}^\dagger:\mathrm{Vec}_{K(s)}\ar@<0.6ex>[r]& \ar@<0.6ex>[l] \mathrm{Isoc}^\dagger\left(X_s/K(s)\right):H^0_\mathrm{rig}(X_s,-)}.
\end{equation}  
on the fibre over $s$. Thus exactly as in the proof of Proposition \ref{trivfibiso}, if $i_s^*E$ is trivial, the counit $g^*g_*E\rightarrow E$ is an isomorphism on the fibre over $s$, and hence an isomorphism. Similarly, since $H^0_\rig(X_s,i_s^*E)\cong \mathrm{Hom}_{\mathrm{Isoc}^\dagger(X/K(s))}(\cur{O}_{X_s},i_s^*E)$, (see Proposition \ref{isoconetwo} below) exactly the same argument as in Proposition \ref{maxsubfib} shows that in general $H^0_\rig(X_s,i_s^*E)\otimes_{K(s)}\cur{O}_{X_s}$ is the largest trivial subobject of $i_s^*E$. Hence if we let $E_0=g^*g_*E$, then $i_s^*E_0\cong H^0_\rig(X_s,i_s^*E)\otimes_{K(s)}\cur{O}_{X_s}$ is the largest trivial sub-object of $i_s^*E$, proving (2), and if $i_s^*E$ is trivial, then $E\cong E_0$, proving (1).
\end{proof}

We now turn our attention to Theorem \ref{two}.

\begin{proposition} \label{isoconetwo} Suppose that $E,E'\in\mathrm{Isoc}^\dagger(X_s/K(s))$. Then there are canonical isomorphisms
\begin{align} \mathrm{Hom}_{\mathrm{Isoc}^\dagger(X_s/K(s))}(E,E') &\cong H^0_\rig(X_s,\cur{H}\mathrm{om}(E,E')) \\
\mathrm{Ext}_{\mathrm{Isoc}^\dagger(X_s/K(s))}(E,E')&\cong H^1_\rig(X_s,\cur{H}\mathrm{om}(E,E')) \nonumber
\end{align}
and moreover if $E,E'$ have Frobenius structures, this induces an isomorphism
\begin{equation}\mathrm{Hom}_{F\text{-}\mathrm{Isoc}^\dagger(X_s/K(s))}(E,E')\cong H^0_\rig(X_s,\cur{H}\mathrm{om}(E,E'))^{\phi=1}
\end{equation}  
as well as a surjection
\begin{equation}\mathrm{Ext}_{F\text{-}\mathrm{Isoc}^\dagger(X_s/K(s))}(E,E')\twoheadrightarrow H^1_\rig(X_s,\cur{H}\mathrm{om}(E,E'))^{\phi=1}
\end{equation}  
\end{proposition}

\begin{proof} The first isomorphism is clear, and this immediately implies the third. The second is Proposition 1.3.1 of \cite{CLS99b}, from which the fourth is then easily deduced.
\end{proof}

We define the $U_n$ inductively as follows. $U_1$ will just be $\cur{O}^\dagger_{X_s}$, and $U_{n+1}$ will be the extension of $U_{n}$ by $\cur{O}^\dagger_{X_s}\otimes_{K(s)} H^1_\mathrm{rig}\left(X_s,U_n^\vee\right)^\vee$ corresponding to the identity under the isomorphisms
\begin{align} \mathrm{Ext}_{\mathrm{Isoc}^\dagger\left(X_s/K(s)\right)}&\left(U_n,\cur{O}^\dagger_{X_s}\otimes_{K(s)} H^1_\mathrm{rig}\left(X_s,U_n^\vee\right)^\vee\right) \\ &
\cong H^1_\mathrm{rig}\left(X_s,U_n^\vee \otimes_{K(s)} H^1_\mathrm{rig}\left(X_s,U_n^\vee\right)^\vee\right)  \nonumber \\
&\cong H^1_\mathrm{rig}\left(X_s,U_n^\vee\right)\otimes_{K(s)} H^1_\mathrm{rig}\left(X_s,U_n^\vee\right)^\vee \nonumber \\
&\cong \mathrm{End}_{K(s)}\left(H^1_\mathrm{rig}\left(X_s,U_n^\vee\right)\right). \nonumber
\end{align}
If we look at the long exact sequence in cohomology associated to the short exact sequence $0 \rightarrow U_n^\vee \rightarrow U_{n+1}^\vee \rightarrow \cur{O}^\dagger_{X_s}\otimes_{K(s)} H^1_\mathrm{rig}\left(X_s,U_n^\vee\right)\rightarrow 0$ we get 
\begin{align} 0\rightarrow H^0_\mathrm{rig}\left(X_s,U_n^\vee\right)&\rightarrow H^0_\mathrm{rig}\left(X_s,U_{n+1}^\vee\right)\rightarrow H^1_\mathrm{rig}\left(X_s,U_n^\vee\right) \\
&
 \overset{\delta}{\rightarrow} H^1_\mathrm{rig}\left(X_s,U_n^\vee\right)\rightarrow H^1_\mathrm{rig}\left(X_s,U_{n+1}^\vee\right). \nonumber
\end{align}

\begin{lemma} \label{connidmor}The connecting homomorphism $\delta$ is the identity.
\end{lemma}

\begin{proof}By dualising, the extension 
\begin{equation} 0\rightarrow U_n^\vee \rightarrow U_{n+1}^\vee\rightarrow \cur{O}^\dagger_{X_s}\otimes_{K(s)} H^1_\mathrm{rig}\left(X_s,U_n^\vee\right)\rightarrow 0
\end{equation}  
corresponds to the identity under the isomorphism
\begin{equation}
\mathrm{Ext}_{\mathrm{Isoc}^\dagger\left(X_s\right/K(s))}\left(\cur{O}^\dagger_{X_s}\otimes_{K(s)} H^1_\mathrm{rig}\left(X_s,U_n^\vee\right),U_n^\vee \right) \cong \mathrm{End}_{K(s)}\left(H^1_\mathrm{rig}\left(X_s,U_n^\vee\right)\right)
\end{equation}  
Now the Lemma follows from the fact that, for an extension  $0\rightarrow E\rightarrow F\rightarrow \cur{O}^\dagger_{X_s}\otimes_{K(s)} V\rightarrow 0$ of a trivial bundle by $E$, the class of the extensions under the isomorphism 
\begin{align} \mathrm{Ext}_{\mathrm{Isoc}^\dagger\left(X_s/K(s)\right)} \left(\cur{O}^\dagger_{X_s}\otimes_K V,E \right)&\cong V^\vee\otimes_{K(s)} H^1_\mathrm{rig}\left(X_s,E\right) \\ &\cong \mathrm{Hom}_{K(s)}\left(V,H^1_\mathrm{rig}\left(X_s,E\right)\right) \nonumber 
\end{align} is just the connecting homomorphism for the long exact sequence
\begin{equation} 0\rightarrow H^0_\mathrm{rig}\left(X_s,E\right)\rightarrow H^0_\mathrm{rig}\left(X_s,F\right)\rightarrow V\rightarrow H^1_\mathrm{rig}\left(X_s,E\right).
\end{equation}  
\end{proof}

In particular, any extension of $U_n$ by a trivial bundle $V\otimes_{K(s)} \cur{O}^\dagger_{X_s}$ is split after pulling back to $U_{n+1}$, and $H^0_\rig\left(X_s,U^\vee_{n+1}\right)\cong H^0_\mathrm{dR}\left(X_s,U^\vee_n\right)$. It then follows by induction that $H^0_\rig\left(X_s,U^\vee_n\right)\cong H^0_\rig(X_s,\cur{O}^\dagger_{X_s})\cong K(s)$ for all $n$. 

\begin{definition} Define the unipotent class of $E\in\cur{N}\mathrm{Isoc}^\dagger\left(X_s/K(s)\right)$ inductively as follows. If $E$ is trivial, then we say $E$ has unipotent class 1. If there exists an extension 
\begin{equation} 0\rightarrow V\otimes_{K(s)} \cur{O}^\dagger_{X_s} \rightarrow E\rightarrow E'\rightarrow 0
\end{equation}  
with $E'$ of unipotent class $\leq m-1$, then we say that $E$ has unipotent class $\leq m$. 
\end{definition}

Now let $x=p(s)$, $u_1=1\in x^*\left(U_1\right)=K(s)$, and choose a compatible system of elements $u_n\in x^*\left(U_n\right)$ mapping to $u_1$.

\begin{proposition} Let $F\in\cur{N}\mathrm{Isoc}^\dagger\left(X_s/K(s)\right)$ be an object of unipotent class $\leq m$. Then for all $n\geq m$ and any $f\in x^*\left(F\right)$ there exists a homomorphism $\alpha:U_n\rightarrow F$ such that $\left(x^*\alpha\right)\left(u_n\right)=f$.
\end{proposition}

\begin{proof} As in the characteristic zero case, we copy the proof of Proposition 2.1.6 of \cite{Had10} and use strong induction on $m$. The case $m=1$ is straightforward. For the inductive step, let $F$ be of unipotent class $m$, and choose an exact sequence
\begin{equation} \label{blahblahblah}0\rightarrow E\overset{\psi}{\rightarrow}F\overset{\phi}{\rightarrow}G\rightarrow 0
\end{equation}  
with $E$ trivial and $G$ of unipotent class $<m$. By induction there exists a unique morphism $\beta:U_{n-1}\rightarrow G$ such that $\left(x^*\phi\right)\left(f\right)=\left( x^*\beta\right)\left(u_{n-1}\right)$. Pulling back the extension (\ref{blahblahblah}) first by the morphism $\beta$ and then by the natural surjection $U_n\rightarrow U_{n-1}$ gives an extension of $U_n$ by $E$, which must split, as observed above.
\begin{equation}\xymatrix{
0 \ar[r]&E \ar[r]\ar@{=}[d]& F'' \ar[r]\ar[d]& U_n \ar[r]\ar[d]\ar@/_1pc/[l]& 0 \\
0 \ar[r]&E \ar[r]\ar@{=}[d]&F' \ar[r]\ar[d]& U_{n-1} \ar[r]\ar[d]& 0 \\
0 \ar[r]& E \ar[r]&F \ar[r]& G\ar[r]& 0
}
\end{equation}  
Let $\gamma:U_n\rightarrow F$ denote the induced morphism, then $\left(x^*\phi\right)\left(\left( x^*\gamma \right)\left(u_n\right)-f\right)=0$. Hence there exists some $e\in x^*E$ such that $\left( x^*\psi\right)\left(e\right)=\left(x^*\gamma\right)\left(u_n\right)-f$. Again by induction we can choose $\gamma':U_n\rightarrow E$ with $\left(x^*\gamma'\right)\left(u_n\right)=e$. Finally let $\alpha=\gamma - \psi\circ\gamma'$, it is easily seen that $\left(x^*\alpha\right)(u_n)=f$. 
\end{proof}

\begin{corollary} Every $E$ in $\cur{N}\mathrm{Isoc}^\dagger\left(X_s/K(s)\right)$ is a quotient of $U_n^{\oplus m}$ for some $n,m\in \N$.
\end{corollary}

Recall that we have the higher direct images $\mathbf{R}^ig_*(E)$ for any $E\in F\text{-}\mathrm{Isoc}^\dagger(X/K)$. Thanks to 2.1.4 of  \cite{Car04}, and the compatibilities already noted between tensor products and pull-backs of arithmetic $\cur{D}$-modules and their counterparts for overconvergent $F$-isocrystals, these satisfy a projection formula
\begin{equation} \mathbf{R}^ig_*(E\otimes g^*E')\cong \mathbf{R}^ig_*(E)\otimes E'
\end{equation}
for any $E\in F\text{-}\mathrm{Isoc}^\dagger(X/K)$ and $E'\in F\text{-}\mathrm{Isoc}^\dagger(S/K)$.

If we let $h$ denote the structure morphism of $S$, then the fact that ${h}_+\circ{g}_+=(h\circ g)_+$ implies that there is a Leray spectral sequence relating $\mathbf{R}^i{h}_*$, $\mathbf{R}^j{g}_*$ and $\mathbf{R}^{i+j}(h\circ g)_*$. Since $S$ is an affine curve and hence $H^2_\rig(S,{g}_*E)=0$,  the exact sequence of low degree terms of this spectral sequence reads
\begin{equation} 0\rightarrow H^1_\rig(S,{g}_*E)\rightarrow H^1_\rig(X,E)\rightarrow H^0_\rig(S,\mathbf{R}^1{g}_*E)\rightarrow 0.
\end{equation}  
We are now in a position to inductively extend the $U_n$ to $X$. Let $W_1=\cur{O}^\dagger_{X}$.

\begin{theorem} \label{fininj} There exists an extension $W_{n+1}$ of $W_n$ by $g^*(\mathbf{R}^1{g}_*W_n^\vee)^\vee$ in the category $\cur{N}_{g}F\text{-}\mathrm{Isoc}^\dagger(X/K)$ such that $i_s^*W_{n+1}=U_{n+1}$ and ${g}_*W_{n+1}^\vee\cong \cur{O}_{S}^\dagger$.
\end{theorem}

\begin{proof} The statement and its proof are by induction on $n$, and in order to prove it we strengthen the induction hypothesis by also requiring that there exists a morphism $p^*W_n^\vee\rightarrow \cur{O}_S^\dagger$ such that the composite morphism $\cur{O}_S^\dagger\cong g_*W_n^\vee\cong p^*g^*g_*W_n^\vee\rightarrow p^*W_n^\vee\rightarrow \cur{O}_S^\dagger$ is an isomorphism.

To check the base case we simply need to verify that ${g}_*\cur{O}_{X}^\dagger\cong \cur{O}_{S}^\dagger$. By the results of the previous section, we get a natural morphism $\cur{O}_S^\dagger\rightarrow g_*\cur{O}^\dagger_X$ as the unit of the adjunction between $g_*$ and $g^*$. By naturality, restricting this morphism to the fibre over $s$ gives us the unit $K(s)\rightarrow H^0_\rig(X_s,\cur{O}_{X_s}^\dagger)$ of the adjunction between $H^0_\rig(X_s,\cdot)$ and $\cdot \otimes_K \cur{O}_{X_s}^\dagger$, which is easily checked to be an isomorphism. Hence by rigidity, $\cur{O}_S^\dagger\rightarrow g_*\cur{O}_X^\dagger$ is an isomorphism.

So now suppose that we have $W_n$ as claimed. We look at the extension group
\begin{equation} \mathrm{Ext}_{F\text{-}\mathrm{Isoc}^\dagger(X/K)}(W_n,g^*(\mathbf{R}^1{g}_*W_n^\vee)^\vee)\twoheadrightarrow H^1_\rig(X, W_n^\vee\otimes_{\cur{O}^\dagger_{X}}g^*(\mathbf{R}^1{g}_*W_n^\vee)^\vee)^{\phi=1}.
\end{equation}  
The Leray spectral sequence, the projection formula above and the induction hypothesis that ${g}_*W_n^\vee\cong \cur{O}_{S}^\dagger$ give us a short exact sequence
\begin{align}\label{hard} 0\rightarrow H^1_\rig(S,(\mathbf{R}^1{g}_*W_n^\vee)^\vee)\rightarrow H^1_\rig(X, W_n^\vee\otimes_{\cur{O}^\dagger_{X}}g^*(\mathbf{R}^1{g}_*W_n^\vee)^\vee)  \\ \rightarrow H^0_\rig(S,\cur{E}\mathrm{nd}(\mathbf{R}^1{g}_*W_n^\vee))\rightarrow 0 \nonumber
\end{align}
which we claim splits compatibly with Frobenius actions. Indeed, pulling back to $S$ via $p$ gives us a map
\begin{equation} H^1_\rig(X,W_n^\vee\otimes_{\cur{O}^\dagger_{X}}g^*(\mathbf{R}^1{g}_*W_n^\vee)^\vee)\rightarrow H^1_\rig(S,p^*W_n^\vee\otimes_{\cur{O}^\dagger_{S}}(\mathbf{R}^1{g}_*W_n^\vee)^\vee)
\end{equation}  
which is again compatible with Frobenius. The projection $p^*W_n^\vee\rightarrow \cur{O}^\dagger_{S}$ induces a map
\begin{equation}H^1_\rig(X,p^*W_n^\vee\otimes_{\cur{O}^\dagger_{S}}(\mathbf{R}^1{g}_*W_n^\vee)^\vee)\rightarrow H^1_\rig(S,(\mathbf{R}^1{g}_*W_n^\vee)^\vee)
\end{equation}  
which is Frobenius compatible, and is such that the composite (dotted) arrow
\begin{equation} \xymatrix{ H^1_\rig(S,(\mathbf{R}^1{g}_*W_n^\vee)^\vee) \ar@{-->}[d] \ar[r] & H^1_\rig(X,W_n^\vee\otimes_{\cur{O}^\dagger_{X}}g^*(\mathbf{R}^1{g}_*W_n^\vee)^\vee) \ar[d] \\
H^1_\rig(S,(\mathbf{R}^1{g}_*W_n^\vee)^\vee) & 
H^1_\rig(S,p^*W_n^\vee\otimes_{\cur{O}^\dagger_{S}}(\mathbf{R}^1{g}_*W_n^\vee)^\vee) \ar[l]}
\end{equation}  
is an isomorphism. Indeed, once the $H^1$'s have been identified with extension groups, the dotted arrow corresponds to push-out along the composite arrow $\cur{O}_S^\dagger\cong g_*W_n^\vee\cong p^*g^*g_*W_n^\vee\rightarrow p^*W_n^\vee\rightarrow \cur{O}_S^\dagger$, which is an isomorphism by the induction hypothesis. Thus the sequence (\ref{hard}) splits as claimed. Let 
\begin{equation}V\subset H^1_\rig(X,W_n^\vee\otimes_{\cur{O}^\dagger_{X}}g^*(\mathbf{R}^1{g}_*W_n^\vee)^\vee)\end{equation}  
be a complementary subspace to $H^1_\rig(S,(\mathbf{R}^1{g}_*W_n^\vee)^\vee)$. By naturality of the Leray spectral sequence we have a commutative diagram
\begin{equation} \xymatrix{ V \ar[d]\ar[r] &H^0_\rig(S,\cur{E}\mathrm{nd}(\mathbf{R}^1{g}_*W_n^\vee)) \ar[d] \\
H^1_\rig(X_s,U_n^\vee\otimes_{K(s)} H^1_\rig(X_s, U_n^\vee)^\vee) \ar@{=}[r] & \mathrm{End}_{K(s)}(H^1_\rig(X_s,U_n^\vee)) }
\end{equation}  
where the left hand vertical arrow is given by restriction to the fibre $X_s$, and the top arrow is an isomorphism. Moreover, all arrows in this diagram are compatible with Frobenius.

The identity in $\mathrm{End}_{K(s)}(H^1_\rig(X_s,U_n^\vee))$, which is Frobenius invariant and corresponds to the extension $U_{n+1}$, lifts to the identity in $H^0_\rig(S,\cur{E}\mathrm{nd}(\mathbf{R}^1{g}_*W_n^\vee))=\mathrm{End}_{\mathrm{Isoc}^\dagger(S)}(\mathbf{R}^1{g}_*W_n^\vee))$, and this element is also Frobenius invariant. Since the upper horizontal map is an isomorphism, compatible with the Frobenius action, we can find a Frobenius invariant class in $V$ mapping to the identity. We let $W_{n+1}'$ be any corresponding extension (the map from the extension group as $F$-isocrystals to the Frobenius invariant part of $H^1$ is surjective). Now, we have a natural map
\begin{equation} \mathrm{Ext}_{F\text{-}\mathrm{Isoc}^\dagger(S/K)}(\cur{O}^\dagger_{S},(\mathbf{R}^1g_*W_n^\vee)^\vee)\overset{g^*}{\rightarrow}  \mathrm{Ext}_{F\text{-}\mathrm{Isoc}^\dagger(X/K)}(W_n,g^*(\mathbf{R}^1g_*W_n^\vee)^\vee)
\end{equation}   
which has a section (denoted $p^*$) induced by the map $p^*W_n^\vee\rightarrow \cur{O}_S^\dagger$, and such that whole diagram
\begin{equation} \xymatrix{ H^1_\rig(S,(\mathbf{R}^1{g}_*W_n^\vee)^\vee) \ar@/_2pc/[r] & H^1_\rig(X,W_n^\vee\otimes_{\cur{O}^\dagger_{X}}g^*(\mathbf{R}^1{g}_*W_n^\vee)^\vee)\ar@/_2pc/[l]\\ \\
 \mathrm{Ext}_{F\text{-}\mathrm{Isoc}^\dagger(S/K)}(\cur{O}^\dagger_{S},(\mathbf{R}^1g_*W_n^\vee)^\vee) \ar[uu]\ar@/_2pc/[r]& \mathrm{Ext}_{F\text{-}\mathrm{Isoc}^\dagger(X/K)}(W_n,g^*(\mathbf{R}^1g_*W_n^\vee)^\vee)\ar[uu]\ar@/_2pc/[l] \\ \\
}
\end{equation}  
commutes. We let $W_{n+1}$ be the extension corresponding to $[W_{n+1}']-g^*p^*[W_{n+1}']$ in  $\mathrm{Ext}_{F\text{-}\mathrm{Isoc}^\dagger(X/K)}(W_n,g^*(\mathbf{R}^1g_*W_n^\vee)^\vee)$. Note that this splits when we pullback via $p^*$ and then push-out via $p^*W_n^\vee\rightarrow \cur{O}_S^\dagger$, and also has the same image as $W_{n+1}'$ inside $H^1_\rig(X,W_n^\vee\otimes_{\cur{O}^\dagger_{X}}g^*(\mathbf{R}^1{g}_*W_n^\vee)^\vee)$.

 To complete the induction we need to show that ${g}_*W_{n+1}^\vee\cong \cur{O}_{S}^\dagger$, and that there exists a map $p_*W_{n+1}^\vee\rightarrow \cur{O}_S^\dagger$ as claimed. We have an exact sequence (using the projection formula and the fact that ${g}_*\cur{O}_{X}^\dagger\cong\cur{O}^\dagger_{S}$)
\begin{equation} 0\rightarrow {g}_*W_n^\vee\rightarrow {g}_*W_{n+1}^\vee\rightarrow \mathbf{R}^1{g}_*W_n^\vee\rightarrow \ldots
\end{equation}  
and it follows from Lemma \ref{connidmor} together with base change that the arrow $g_*W_n^\vee\rightarrow g_*W_{n+1}^\vee$ restricts to an isomorphism on the fibre at $s$. Thus by rigidity it is an isomorphism. Finally, we have an exact sequence
\begin{equation} 0 \rightarrow p^*W_n^\vee\rightarrow p^*W_{n+1}\rightarrow (\mathbf{R}^1g_*W_n^\vee)^\vee\rightarrow 0
\end{equation}   which  splits when we push-out via the map $p^*W_n^\vee\rightarrow \cur{O}_S^\dagger$. This splitting induces a map $p^*W_{n+1}^\vee\rightarrow \cur{O}_S^\dagger$ such that the diagram
\begin{equation}\xymatrix{
p^*W_n^\vee \ar[dr]\ar[r] & p^*W_{n+1}^\vee\ar[d] \\
& \cur{O}_S^\dagger}
\end{equation}  
commutes. Now the fact that the diagram
\begin{equation} \xymatrix{
& g_*W_{n+1}^\vee \ar[r]& p^*g^*g_*W_{n+1}^\vee \ar[r] & p^*W_{n+1}^\vee \ar[dr] \\
\cur{O}_S^\dagger \ar[r]\ar[ur] & g_*W_n^\vee\ar[u] \ar[r]& p^*g^*g_*W_n^\vee \ar[r]\ar[u] & p^*W_n^\vee \ar[r]\ar[u] & \cur{O}_S^\dagger
}
\end{equation}  
commutes implies that the composite along the top row is an isomorphism, finishing the proof.

\end{proof}

To complete the proof of Theorem \ref{two}, we use the base extension functor 
\begin{equation}-\otimes_K K(s):\cur{N}_{g}F\text{-}\mathrm{Isoc}^\dagger(X/K)\rightarrow \cur{N}_{g}F\text{-}\mathrm{Isoc}^\dagger(X/K)\otimes_K K(s),
\end{equation}
which is defined on pages 155-156 of \cite{MD81}, to view the $W_n$ as objects of the latter category.

\subsection{Extension to proper curves, Frobenius structures}\label{frobstrut}

In this section we use gluing methods to define $\pi_1^\rig(X/S,p)$ whenever $S$ is a smooth, geometrically connected curve over $k$. Since we will depend on the results from the previous section, we will assume that Hypothesis \ref{hypothesis} holds Zariski locally on $S$.

\begin{lemma} Let $j:T\rightarrow S$ be a morphism of smooth, geometrically connected affine curves over $k$. Then the canonical morphism $\pi_1^\rig(X_T/T,p_T)\rightarrow j^*(\pi_1^\rig(X/S,p))$ is an isomorphism.
\end{lemma}

\begin{proof} By rigidity, it suffices to show that it is an isomorphism on stalks. But this follows from the fact that the induced map on stalks is just the canonical isomorphism $\pi_1^\rig((X_T)_t,(p_T)_t)\isomto \pi_1^\rig(X_{j(t)},p_{j(t)})$.
\end{proof}

Now suppose that $S$ is a (not necessarily affine) smooth, geometrically connected curve. Let $\{S_i\}$ be a cover of $S$ by affine curves, let $g_i:X_i\rightarrow S_i$ be the pull-back of $g$ to $S_i$, and $p_i$ the induced section. Let $S_{ij}=S_i\cap S_j$, and similarly denote $g_{ij},X_{ij},p_{ij}$. The category $F\text{-}\mathrm{Isoc}^\dagger(S/K)$ is Zariski-local on $S$, and the above lemma shows that we have isomorphisms 
\begin{equation}\pi_1^\rig(X_i/S_i,p_i)|_{S_{ij}}\cong \pi_1^\rig(X_{ij}/S_{ij},p_{ij})\cong \pi_1^\rig(X_j/S_j,p_j)|_{S_{ij}}
\end{equation}  
for all $i,j$. These satisfy the co-cycle condition on triple intersections, and hence glue to give an affine group scheme $\pi_1^\rig(X/S,p)$ over $F\text{-}\mathrm{Isoc}^\dagger(S/K)$. Using the above lemma, it is easy to check that this object is independent of the choice of affine covering $\{S_i\}$, up to canonical isomorphism.

\begin{definition} When $S$ is a smooth, geometrically connected curve, we will denote by $\pi_1^\rig(X/S,p)$ the affine group scheme just constructed by gluing, and not the object defined in previous sections.
\end{definition}

Now let $f:T\rightarrow S$ be a morphism of curves, smooth and geometrically connected over finite fields $k'$ and $k$ respectively. Let $K'$ denote the unique unramified extension of $K$ with residue field $k'$.

\begin{lemma} \label{bsch} \begin{enumerate} \item Let $s\in S$ be a closed point. Then there is an isomorphism $\pi_1^\rig(X/S,p)_s\cong \pi_1^\rig(X_s,p_s)$. 
\item There is a natural isomorphism $\pi_1^\rig(X_T/T,p_T)\isomto f^*(\pi_1^\rig(X/S,p))$.
\end{enumerate}
\end{lemma}

\begin{proof} The first immediately follows from the corresponding result when $S,T$ are affine. The second follows from the first.
\end{proof}

\begin{remark} If $x\in X(S)$ is a(nother) point, then by exactly the same technique we can glue the path torsors over affine sub-curves of $S$ to obtain path torsors under $\pi_1^\rig(X/S,p)$.
\end{remark}

The upshot of the previous section is that we now have an affine group scheme $\pi_1^\rig(X/S,p)$ over the Tannakian category $F\text{-}\mathrm{Isoc}^\dagger(S/K)$ whose fibre (ignoring Frobenius structures) over any closed point $s$ is the usual rigid fundamental group $\pi_1^\rig(X_s,p_s)$ as defined by Chiarellotto and le Stum in \cite{CLS99a}. In Chapter II of \cite{Chi98}, Chiarellotto defines a Frobenius isomorphism $F_*:\pi_1^\rig(X_s,p_s)\isomto\pi_1^\rig(X_s,p_s)$, by using the fact that Frobenius pullback induces an automorphism of the category $\cur{N}\mathrm{Isoc}^\dagger(X_s/K)$. Since we have constructed $\pi_1^\rig(X/S,p)$ as an affine group scheme over $F\text{-}\mathrm{Isoc}^\dagger(S/K)$, it comes with a Frobenius structure that we can compare with Chiarellotto's. However, it is not obvious to us exactly what the relationship between these two Frobenius structures is, so instead we will endow $\pi_1^\rig(X/S,p)$ with a different Frobenius, which we will be able to compare with the natural Frobenius on the fibres.

\begin{remark} From now onward, we will consider $\pi_1^\rig(X/S,p)$ as an affine group scheme over $\mathrm{Isoc}^\dagger(S/K)$, via the forgetful functor. Note that Lemma \ref{bsch} still holds, \emph{a fortiori}, if we ignore the $F$-structure.
\end{remark}

Let $\sigma_S:S\rightarrow S$ denote the $k$-linear Frobenius, $X'=X\times_{S,\sigma_S} S$ the base change of $X$ by $\sigma_S$, and $\sigma_{X/S}:X\rightarrow X'$ the relative Frobenius induced by the $k$-linear Frobenius $\sigma_X$ of $X$. Let $p'$ be the induced point of $X'$, and $q=\sigma_{X/S}\circ p\in X'(S)$. Then by functoriality and base change we get a homomorphism
\begin{equation} \pi_1^\rig(X/S,p)\rightarrow \pi_1^\rig(X'/S,q)
\end{equation}  
and an isomorphism
\begin{equation} \pi_1^\rig(X'/S,p')\isomto \sigma_S^*\pi_1^\rig(X/S,p).
\end{equation} 
One can easily check that $p'=q\in X(S)$, and hence we get a natural morphism  $\phi: \pi_1^\rig(X/S,p)\rightarrow \sigma_S^*\pi_1^\rig(X/S,p)$.

\begin{lemma} This is an isomorphism.
\end{lemma}

\begin{proof} Let $s\in S$ be a closed point, with residue field $k(s)$ of size $q^a$. The map induced by $\phi^a$ on the fibre $\pi_1^\rig(X_s,p_s)$ over $s$ is the same as that induced by pulling back unipotent isocrystals on $X_s$ by the $k(s)$-linear Frobenius on $X_s$. This is proved in Chapter II of \cite{Chi98} to be an isomorphism, thus $\phi^a$ is an isomorphism by rigidity. Hence $\phi$ is also an isomorphism.
\end{proof}

We now let $F_*:\sigma_S^*\pi_1^\rig(X/S,p)\isomto\pi_1^\rig(X/S,p)$ denote the inverse of $\phi$, which by the proof of the previous lemma, reduces to the Frobenius structure as defined by Chiarellotto on closed fibres.

\begin{definition} When we refer to `the' Frobenius on $\pi_1^\rig(X/S,p)$, we will mean the isomorphism $F_*$ just defined.
\end{definition}

\subsection{Cohomology and period maps}\label{cohomperiod}

In this section we study the non-abelian cohomology of the unipotent quotients $\pi_1^\rig(X/S,p)_n$ of $\pi_1^\rig(X/S,p)$. Assumptions and notations will be exactly as in the previous two sections. Recall from Section \ref{reltanna} the notion of a torsor under an affine group scheme $U$ over $\mathrm{Isoc}^\dagger(S/K)$.

\begin{definition} We define $H^1_\rig(S,U)$ to be the pointed set of isomorphism classes of torsors under $U$.
\end{definition}

\begin{example} \label{vectorgroup} Suppose that $U$ is the vector scheme associated to an overconvergent isocrystal $E$. Then Exemple 5.10 of \cite{Del89} shows that there is a bijection $H^1_\rig(S,U)\isomto H^1_\rig(S,E)$.
\end{example}

If $U$ has a Frobenius structure, that is an isomorphism $\phi:\sigma_S^*U\isomto U$, then we can define an $F$-torsor under $U$ to be a $U$-torsor $P$, together with a Frobenius isomorphism $\phi_P:\sigma_S^*P\isomto P$ such that the action map $P\times U\rightarrow P$ is compatible with Frobenius.

\begin{definition}  We define $H^1_{F,\rig}(S,U)$ to be the set of isomorphism classes of $F$-torsors under $U$.
\end{definition}

Given any torsor $P$ under $U$, without $F$-structure, $\sigma_S^*P$ will be a torsor under $\sigma_S^*U$, and hence we can use the isomorphism $\phi$ to consider $\sigma_S^*P$ as a torsor under $U$. Hence we get a Frobenius action $\phi:H^1_\rig(S,U)\rightarrow H^1_\rig(S,U)$, and it is easy to see that the forgetful map
\begin{equation} H^1_{F,\rig}(S,U)\rightarrow H^1_\rig(S,U)
\end{equation}  
is a surjection onto the subset $H^1_\rig(S,U)^{\phi=\mathrm{id}}$ fixed by the action of $\phi$.

Given any point $x\in X(S)$, we have the path torsors $P(x)$ under $\pi_1^\rig(X/S,p)$ as well as the finite level versions $P(x)_n$. Moreover, these come with Frobenius structures, and hence we get compatible maps
\begin{equation} \xymatrix{ X(S)\ar[r]\ar[dr] & H^1_{F,\rig}(S,\pi_1^\rig(X/S,p)_n)\ar[d] \\
& H^1_{\rig}(S,\pi_1^\rig(X/S,p)_n)^{\phi=\mathrm{id}}}
\end{equation}  
for each $n\geq1$. 

In order to get a handle on this `non-abelian' $H^1$, we first discuss the generalisation of Theorem 2.11 of \cite{MD81} to non-neutral Tannakian categories via groupoids and their representations, following \cite{Del90}. The reason for doing this is to obtain a generalisation of Example \ref{actiontors} giving a more explicit description of $H^1_\rig(S,U)$. 

So let $K$ be a field of characteristic $0$, and $Y$ a $K$-scheme.

\begin{definition} A $K$-groupoid acting on $Y$ is a $K$-scheme $G$, together with `source' and `target' morphisms $s,t:G\rightarrow Y$ and a `law of composition' $\circ : G\times_{^s Y^t} G\rightarrow G$, which is a morphism of $Y\times_K Y$-schemes ($G\times_{^sY^t}G$ considered as a $Y\times_K Y$ scheme via the composition of the projection to $S$ with the diagonal $Y\rightarrow Y\times _K Y$, $G$ considered as a $Y\times_K Y$-scheme via $s\times t$) such that the following conditions hold. For any $K$-scheme $T$, the data of $Y(T)$, $G(T)$, $s,t,\circ$ forms a groupoid, where $Y(T)$ is the set of objects and $G(T)$ the set of morphisms.
\end{definition}

\begin{example} Suppose that $Y=\spec{K}$. Then a $K$-groupoid acting on $Y$ is just a group scheme over $K$.
\end{example}

\begin{definition}\label{groupoidrep} If $G$ is a $K$-groupoid acting on $Y$, then a representation of $G$ is a quasi-coherent $\cur{O}_Y$-module $V$, together with a morphism $\rho(g):s(g)^*V\rightarrow t(g)^*V$ for any $K$-scheme $T$ and any point $g\in G(T)$. These morphisms must be compatible with base change $T'\rightarrow T$, as well as with the law of composition on $G$. Finally, if $\mathrm{id}_y\in G(T)$ is the `identity morphism' corresponding to the `object' $y\in Y(T)$, then we require the morphism $\rho(\mathrm{id}_y)$ to be the identity. A morphism of representations is defined in the obvious way, and we denote the category of \emph{coherent} representations by $\mathrm{Rep}(Y:G)$. Of course we can similarly define actions of $G$ on any (group) scheme $U$ over $Y$, by instead requiring morphisms $\rho(g):U\times_{Y,s(g)}T\rightarrow U\times_{Y,t(g)}T$ of (group) schemes over $T$.
\end{definition}

\begin{example} If $Y=\spec{K}$, then this just boils down to the usual definition of a representation of a group scheme over $K$.
\end{example}

Now suppose that $\cur{C}$ is a Tannakian category over $K$, which admits a fibre functor $\omega:\cur{C}\rightarrow \mathrm{Vec}_L$ taking values in some finite extension $L/K$. Let $\mathrm{pr}_i:\spec{L\otimes_K L}\rightarrow \spec{L}$ for $i=1,2$ denote the two projections. Then we get two fibre functors $\mathrm{pr}_i^*\circ\omega:\cur{C}\rightarrow \mathrm{Mod}_{\mathrm{f.g.}}(L\otimes_K L)$ taking values in the category of finitely generated $L\otimes_K L$-modules, and the functor of isomorphisms $\underline{\mathrm{Isom}}^\otimes(\mathrm{pr}_1^*\circ\omega,\mathrm{pr}_2^*\circ\omega)$ is represented by an affine scheme $\underline{\mathrm{Aut}}^\otimes_K(\omega)$ over $L\otimes_K L$. The composite of the map $\underline{\mathrm{Aut}}^\otimes_K(\omega)\rightarrow \spec{L\otimes_KL}$ with the two projections to $\spec{L}$ makes $\underline{\mathrm{Aut}}^\otimes_K(\omega)$ into a $K$-groupoid acting on $\spec{L}$. Moreover, if $E$ is an object of $\cur{C}$, then $\omega(E)$ becomes a representation of $\underline{\mathrm{Aut}}^\otimes_K(\omega)$ in the obvious way. Thus we get a functor
\begin{equation} \cur{C}\rightarrow \mathrm{Rep}(L:\underline{\mathrm{Aut}}^\otimes_K(\omega))
\end{equation}  
and Th\'{e}or\`{e}me (1.12) of \cite{Del90} states (in particular) the following.

\begin{theorem} The induced functor $\cur{C}\rightarrow \mathrm{Rep}(L:\underline{\mathrm{Aut}}^\otimes_K(\omega))$ is an equivalence of Tannakian categories.
\end{theorem}

Finally, to get the generalisation of Example \ref{actiontors} that we need, the following technical lemma is necessary.

\begin{lemma} (\cite{Del90}, Corollaire 3.9). Let $L/K$ be finite, and $G$ a $K$-groupoid acting on $\spec{L}$, affine and faithfully flat over over $L\otimes_K L$. Then any representation $V$ of $G$ is the colimit of its finite dimensional sub-representations.
\end{lemma}

\begin{corollary} If $\cur{C}$ is a Tannakian category over $K$, $\omega$ a fibre functor with values in $L$, then an affine (group) scheme over $\cur{C}$ `is' just an affine (group) scheme over $L$ together with an action of $\underline{\mathrm{Aut}}^\otimes_K(\omega)$, and morphism of such objects `are' just $\underline{\mathrm{Aut}}^\otimes_K(\omega)$-equivariant morphisms.
\end{corollary}

\begin{definition}
Let $G$ be a $K$-groupoid acting on $\spec{L}$. If $U$ is a group scheme over $L$ with a $G$-action, we will denote by $H^1(G,U)$ the set of isomorphism classes of $G$-equivariant torsors under $U$.
\end{definition}

\begin{example}
\begin{itemize}
\item If $V$ is a representation of $G$, then $\spec{\mathrm{Sym}(V^\vee)}$ naturally becomes a group scheme over $L$ with a $G$-action. We will refer to this latter object as the vector scheme associated to $V$.
\item If $U$ is a unipotent affine group scheme over $\mathrm{Isoc}^\dagger(S/K)$ as above, then for any closed point $s\in S$, the unipotent group $U_s$ over $K(s)$ attains an action of the $K$-groupoid $\underline{\mathrm{Aut}}^\otimes_K(s^*)$, and there is a natural bijection of sets
\begin{equation} H^1_\rig(S,U)\isomto H^1(\underline{\mathrm{Aut}}^\otimes_K(s^*),U_s).
\end{equation}  
\end{itemize}
\end{example}

Suppose that $Y=\spec{L}$, with $L/K$ finite, and let $G$ be a $K$ groupoid acting on $Y$. Let $U$ be a unipotent group over $L$, on which $G$ acts.

\begin{definition} A 1-cocyle for $G$ with values in $U$ is a map of $K$-schemes $\phi:G\rightarrow U$ such that
\begin{itemize} \item The diagram
\begin{equation}
\xymatrix{
G\ar[dr]_t\ar[r]^\phi & U \ar[d]^{\mathrm{canonical}} \\
& \spec{L}
} 
\end{equation}   commutes.
\item For any $K$-scheme $T$, and points $g,h\in G(T)$ which are composable in the sense that $s(g)=t(h)$, $\phi(gh)=\phi(g)\cdot\rho(g)(\phi(h))$ holds. This equality needs some explaination. By the first condition above, $\phi(g)$ lands in the subset $\mathrm{Hom}_T(T,U\times_{L,t(g)}T)$ of $\mathrm{Hom}_K(T,U)$ which consists of those morphisms $T\rightarrow U$ which are such that the diagram
\begin{equation}
\xymatrix{
T\ar[dr]_{t(g)}\ar[r] & U \ar[d]^{\mathrm{canonical}} \\
& \spec{L}
} 
\end{equation}
commutes. Similarly, $\phi(h)\in \mathrm{Hom}_T(T,U\times_{L,t(h)}T)=\mathrm{Hom}_T(T,U\times_{L,s(g)}T)$. Since $U/L$ is a group scheme, $\mathrm{Hom}_T(T,U\times_{L,t(g)}T)$ is a group, and the action of $G$ on $U$ gives a homomorphism 
\begin{equation}
\rho(g):\mathrm{Hom}_T(T,U\times_{L,s(g)}T)\rightarrow \mathrm{Hom}_T(T,U\times_{L,t(g)}T).
\end{equation}
Hence the equality $\phi(gh)=\phi(g)\cdot\rho(g)(\phi(h))$ makes sense inside the group  $\mathrm{Hom}_T(T,U\times_{L,t(g)}T)$.
\end{itemize}  
The set of 1-cocycles with coefficients in $U$ is denoted $Z^1(G,U)$. This set has a natural action of $U(L)$ via 
\begin{equation}
(\phi* u)(g)=(t(g)^*u)^{-1}\cdot\phi(g)\cdot \rho(g)(s(g)^*(u))
\end{equation} for any $g\in G(T)$, as above this makes sense inside the group $\mathrm{Hom}_T(T,U\times_{L,t(g)}T)$.
\end{definition}

The point of introducing these definitions is the following.

\begin{lemma} There is a bijection between the non-abelian cohomology set $H^1(G,U)$ and the set of orbits of $Z^1(G,U)$ under the action of $U(L)$.
\end{lemma}

\begin{proof} Let $P$ be a $G$-equivariant torsor under $U$. Since any torsor under a unipotent group scheme over an affine scheme is trivial, we may choose a point $p\in P(L)$. Now, for any $g\in G(T)$ we can consider the points $t(g)^*p$ and $s(g)^*p$ inside $\mathrm{Hom}_T(T,P\times_{L,t(g)}T)$ and $\mathrm{Hom}_T(T,P\times_{L,s(g)}T)$ respectively. We have a morphism $\rho(g):P\times_{L,s(g)}T\rightarrow P\times_{L,t(g)}T$ and hence there exists a unique element $\phi(g)\in U\times_{L,t(g)}T (T)$ such that $t(g)^*p\phi(g)=\rho(g)s(g)^*p$. Thus we get some $\phi(g)\in U(T)$, and the map $g\mapsto\phi(g)$ is functorial, giving a map of schemes $\phi:G\rightarrow U$. The fact that $\phi(g)\in \mathrm{Hom}_T(T,U\times_{L,t(g)}T)$ means that the diagram
\begin{equation}
\xymatrix{
G\ar[dr]_t\ar[r]^\phi & U \ar[d]^{\mathrm{canonical}} \\
& \spec{L}
} 
\end{equation}   commutes, and one easily checks that $\phi$ satisfies the cocycle condition. A different choice of $p$ differs by an element of $U(L)$, and one easily sees that this modifies $\phi$ exactly as in the action of $U(L)$ on $Z^1(G,U)$ defined above. Hence we get a well defined map
\begin{equation} H^1(G,U)\rightarrow Z^1(G,U)/U(L).
\end{equation}  
Conversely, given a cocycle $\phi:G\rightarrow U$, we can define a torsor $P$ as follows. The underlying scheme of $P$ is just $U$, and the action of $U$ on $P$ is just the usual action of right multiplication. We use the cocycle $\phi$ to twist the action of $G$ as follows. If $g\in G(T)$, then we define $\rho(g):P\times_{L,s(g)}T\rightarrow P\times_{L,t(g)}T$ to be the unique map, compatible with the $U$ action, taking the identity of $U\times_{L,s(g)}T=P\times_{L,s(g)}T$ to $\phi(g)\in U\times_{L,t(g)}T=P\times_{L,t(g)}T$. One easily checks that this descends to the quotient $Z^1(G,U)/U(L)$, and provides an inverse to the map defined above.
\end{proof}

We now want to investigate more closely the case when $U$ is a vector scheme, coming from some finite dimensional representation $V$ of $G$. In this case we define, for any $n\geq0$ the space $C^n(G,V)$ of $n$-cochains of $G$ in $V$ as follows. Let $G^{(n)}$ denote the scheme of `$n$-fold composable arrows in $G$', that is the sub-scheme of $G\times_K \ldots\times_K G$ ($n$ copies), consisting of those points $(g_1,\ldots, g_n)$ such that $s(g_i)=t(g_{i+1})$ for all $i$, by convention we set $G^{(0)}=\spec{L}$. Then the space of $n$-cochains is simply the space of global sections of the coherent sheaf $(\delta_1^n)^*V$ on $G^{(n)}$, where $\delta_1^n:G^{(n)}\rightarrow \spec{L}$ is defined to be  the map $t\circ\mathrm{pr}_1$, where $\mathrm{pr}_1:G^{n}\rightarrow G$ is projection onto the first factor. This can also be viewed as the set of morphisms
$G^{(n)}\rightarrow \spec{\mathrm{Sym}(V^\vee)}$ making the diagram
\begin{equation}
\xymatrix{
G^{(n)}\ar[dr]_{t\circ\mathrm{pr}_1}\ar[r] & \spec{\mathrm{Sym}(V^\vee)} \ar[d]^{\mathrm{canonical}} \\
& \spec{L}
} 
\end{equation}   commute, and hence we can define differentials $d^n:C^n(G,V)\rightarrow C^{n+1}(G,V)$ by 
\begin{align} (d^nf)(g_1,\ldots,g_{n+1})=\rho(g_1)f(g_2,&\ldots,g_{n+1})\\+\sum_{i=1^n}&(-1)^if(g_1,\ldots,g_ig_{i+1},\ldots,g_{n+1}\nonumber)\\&+(-1)^{n+1}f(g_1,\ldots,g_n)\nonumber
\end{align}
for $n\geq1$, where $g_1,\ldots,g_{n+1}$ are composable elements of $G(T)$, and all the summands on the RHS are global sections of the coherent sheaf $t(g_1)^*V$ on $T$. For $n=0$ we define $(d^0f)(g)=\rho(g)f(s(g)-f(t(g))$. It is easily checked that these differentials make $C^\bullet(G,V)$ into a chain complex, and we define the cohomology of $G$ with coefficients in $V$ to be the cohomology of this complex:
\begin{equation} H^n(G,V):=H^n(C^\bullet(G,V)).
\end{equation}  

\begin{lemma} Let $V$ be a representation of the groupoid $G$ acting on $\spec{L}$. Then there is a canonical bijection $H^1(G,V)\isomto H^1(G,\spec{\mathrm{Sym}(V^\vee)})$
\end{lemma}

\begin{proof} Taking into account the description of the latter in terms of cocyles modulo the action of $V$, this is straightforward algebra.
\end{proof}

So far we have been working over a field $K$, however, exactly the same definitions make sense over any $K$-algebra $R$, and we can define the cohomology of an $R$-groupoid acting on $\spec{R\times_K L}$. There is an obvious base extension functor, taking $K$-groupoids to $R$-groupoids, and hence we can define cohomology functors $\underline{H}^n(G,V)$ for any representation $V$ of $G$.

\begin{proposition} Suppose that $G=\spec{A}$ is affine. Then for any $K$-algebra there are a canonical isomorphisms $H^n(G_R,V_R)\isomto H^n(G,V)\otimes_K R$ for all $n\geq0$.
\end{proposition}

\begin{proof} In this case, there is an alternative algebraic description of the complex $C^\bullet(G,V)$. First of all, $A$ is a commutative $L\otimes_K L$-algebra, hence $A$ becomes an $L$-module in two different ways, using the two maps $L\rightarrow L\otimes_K L$. We will refer to these as the `left' and `right' structures, these two different $L$-module structures induce the same $K$-module structure. The groupoid structure corresponds to a morphism $\Delta:A\rightarrow A\otimes_LA$, using the two different $L$-module structures to form the tensor product.

The action of $G$ on  a representation $V$ can be described by an $L$-linear map $\Delta_V:V\rightarrow V\otimes_{L,t} A$, where on the RHS we use the `left' $L$-module structure on $A$ to form the tensor product, and define the $L$-module structure on the result via the `right' $L$-module structure on $A$. This map is required to satisfy axioms analogous to the comodule axioms for the description of a representation of an affine group scheme.

Hence the group $C^n(G,V)$ of $n$-cochains is simply the $L$-module $V\otimes_L A\otimes_L \ldots \otimes_L A$ ($n$ copies of $A$). We can describe the boundary maps $d^n$ algebraically as well by
\begin{align} d^n(v\otimes a_1\otimes\ldots\otimes a_n)=\Delta_V(v)&\otimes a_1\otimes\ldots\otimes a_n \\
+\sum_{i=1}^n& v\otimes a_1\otimes\ldots\otimes\Delta(a_i)\otimes\ldots\otimes a_n \nonumber\\ &+v\otimes a_1\otimes\ldots\otimes a_n\otimes 1. \nonumber
\end{align}
Exactly the same discussion applies over any $K$-algebra $R$, and one immediately sees that there is an isomorphism of complexes $C^\bullet(G_R,V_R)\cong C^\bullet(G,V)\otimes_K R$. Since any $K$-algebra is flat, the result follows.
\end{proof}

\begin{remark} In other words, the cohomology functor $\underline{H}^n(G,V)$ is represented by the vector scheme associated to $H^n(G,V)$.
\end{remark}

If $U$ is a unipotent group scheme on which $G$ acts, we can also extend the set $H^1(G,U)$ to a functor of $K$-algebras in the same way. We can also define $H^0(G,U)$ to be the group of all $u\in U(L)$ such that $\rho(g)s(g)^*u=t(g)^*$ for any $g\in G(T)$, and any $K$-scheme $T$. This also extends to a functor of $K$-algebras in the obvious way. It is straightforward to check that $\underline{H}^0(G,\spec{\mathrm{Sym}(V^\vee)})=\underline{H}^0(G,V)$ whenever $V$ is a representation of $G$.

Recall that if $U$ is a unipotent group scheme, we define $U^n$ inductively by $U^1=[U,U]$ and $U^n=[U^{n-1},U]$ and $U_n$ by $U_n=U/U^n$. Since $U$ is unipotent over $K$, a field of characteristic zero we know that each $U^n/U^{n+1}$ is a vector scheme, and that $U=U_N$ for large enough $N$.

\begin{theorem} Let $U$ be a unipotent group scheme acted on by $G$. Assume that $G$ is affine, and for all $n\geq 1$, $H^0(G,U^n/U^{n+1})=0$. Then the functor $\underline{H}^1(G,U)$ is represented by an affine scheme over $K$.
\end{theorem}

\begin{proof} The hypotheses imply that $\underline{H}^0(G,U^n/U^{n+1})(R)=0$ for all $K$-algebras $R$, and $\underline{H}^0(G,U)(R)=0$.

We will prove the theorem by induction on the unipotence degree of $U$, and our argument is almost word for word that given by Kim in the proof of Proposition 2, Section 1 of \cite{Kim05}. When $U$ is just a vector scheme associated to a representation of $G$, then we already know that $\underline{H}^n(G,U)$ is representable for all $n$. For general $U$, we know that we can find an exact sequence
\begin{equation} 1\rightarrow V\rightarrow U\rightarrow W\rightarrow 1
\end{equation}  
realising $U$ as a central extension of a unipotent group $W$ of lower unipotence degree by a vector scheme $V$. Looking at the long exact sequence in cohomology associated to this exact sequence, the boundary map $H^1(G_R,W_R)\rightarrow H^2(G_R,V_R)$ is a functorial map between representables (using the induction hypothesis for representability of $\underline{H}^1(G,W_R))$ and hence the pre-image of $0\in \underline{H}^2(G,V)$ is an (affine) closed sub-scheme of $\underline{H}^1(G,W)$, which we will denote by $I(G,W)$. Thus we get a vector scheme $\underline{H}^1(G,V)$, an affine scheme $I(G,W)$, and an exact sequence
\begin{equation} 1\rightarrow \underline{H}^1(G,V)(R)\rightarrow \underline{H}^1(G,U)(R)\rightarrow I(G,W)(R)\rightarrow 1
\end{equation}  
for all $R$. We now proceed \emph{exactly} as in the proof of Proposition 2, Section 1 of \cite{Kim05} to obtain an isomorphism of functors $\underline{H}^1(G,U)\cong \underline{H}^1(G,V)\times I(G,W)$, showing that $\underline{H}^1(G,U)$ is an affine scheme over $K$.
\end{proof}

\begin{corollary} With the assumptions as in the previous theorem, assume further that  $H^1(G,U^i/U^{i+1})$ is finite dimensional for each $n$. Then $\underline{H}^1(G,U_n)$ is of finite type over $K$, of dimension at most $\sum_{i=1}^{n-1}\dim_KH^1(G,U^i/U^{i+1})$
\end{corollary}

Recall that for a `good' morphism $f:X\rightarrow S$ over a finite field, with $S$ a curve and $X$ satisfying Hypothesis \ref{hypothesis}, we have the period map
\begin{equation} X(S)\rightarrow H^1_\rig(S,\pi_1(X/S,p)_n)
\end{equation}  
taking a section to the corresponding path torsor. Choosing a closed point $s\in S$ means we can interpret this map as
\begin{equation} X(S)\rightarrow H^1(\underline{\mathrm{Aut}}_K^\otimes(s^*),\pi_1^\rig(X_s,p(s))_n).
\end{equation}  
This latter set has the structure of an algebraic variety over $K$ under the condition that
\begin{equation} H^0_\rig(S,\pi_1^\rig(X/S,p)^n/\pi_1^\rig(X/S,p)^{n+1})
\end{equation}  
is zero for each $n$. If, for example, $X$ is a model for a smooth projective curve $C$ over a function field, then we expect this condition to be satisfied under certain non-isotriviality assumptions on the Jacobian of $C$. 

\section*{Acknowledgements}

This paper is intended to form part of the author's PhD thesis, and was written under the supervision of Ambrus P\'{a}l at Imperial College, London, funded by an EPSRC Doctoral Training Grant. He would like to thank Dr. P\'{a}l both for suggesting the direction of research and for providing support, encouragement and many hours of fruitful discussions, without which this paper would not have been written.

\bibliographystyle{mysty}
\bibliography{/Users/Chris/Dropbox/LaTeX/lib.bib}

\providecommand{\bysame}{\leavevmode\hbox to3em{\hrulefill}\thinspace}
\providecommand{\MR}{\relax\ifhmode\unskip\space\fi MR }
\providecommand{\MRhref}[2]{%
  \href{http://www.ams.org/mathscinet-getitem?mr=#1}{#2}
}
\providecommand{\href}[2]{#2}
\begin{thebibliography}{DMSS00}

\bibitem[Ber96a]{Ber96b}
P.~Berthelot, \emph{Cohomologie rigide et cohomologie ridige {\`a} supports
  propres, premi{\`e}re partie}, preprint (1996).

\bibitem[Ber96b]{Ber96a}
\bysame, \emph{$\cur{D}$-modules arithm{\'e}tique {I}. {O}p{\'e}rateurs
  diff{\'e}rentiels de niveau fini}, Ann. Sci. Ecole. Norm. Sup. \textbf{29}
  (1996), 185--272.

\bibitem[Ber00]{Ber00}
\bysame, \emph{$\cur{D}$-modules arithm{\'e}tique {II}. {D}escente par
  {F}robenius}, M{\'e}moires de la SMF (2000), no.~81.

\bibitem[Ber02]{Ber02}
\bysame, \emph{Introduction {\`a} la th{\'e}orie arithm{\'e}tique des
  {$\mathscr{D}$}-modules}, Cohomologies {$p$}-adiques et applications
  arithm{\'e}tiques, II, no. 279, Asterisque, 2002, pp.~1--80.

\bibitem[Bor87]{Bor87}
A.~Borel, \emph{Algebraic {D}-modules}, Perspectives in Mathematics, vol.~2,
  Academic Press Inc., 1987.

\bibitem[Car04]{Car04}
D.~Caro, \emph{$\cur{D}$-modules arithm\'{e}tiques surcoh\'{e}rents.
  {A}pplications aux fonctions {L}}, Ann. Inst. Fourier, Grenoble \textbf{54}
  (2004), no.~6, 1943--1996.

\bibitem[Car06]{Car06a}
\bysame, \emph{D\'{e}vissages des {$F$}-complexes de {$\cur{D}$}-modules
  arithm\'{e}tiques en {$F$}-isocristaux surconvergents}, Inventionnes
  Mathematica \textbf{166} (2006), 397--456.

\bibitem[Car07]{Car07}
\bysame, \emph{{$F$}-isocristaux surconvergents et surcoh\'{e}rence
  differentielle}, Inventionnes Mathematica \textbf{170} (2007), 507--539.

\bibitem[Car09]{Car09b}
\bysame, \emph{$\cur{D}$-modules arithm\'{e}tiques surholonomes}, Ann. Sci.
  Ecole. Norm. Sup. \textbf{42} (2009), no.~1, 141--192.

\bibitem[Car15a]{Car15a}
\bysame, \emph{Sur la pr\'eservation de la surconvergence par l'image directe
  d'un morphisme propre et lisse}, Ann. Sci. \'Ec. Norm. Sup\'er. (4)
  \textbf{48} (2015), no.~1, 131--169.

\bibitem[Car15b]{Car15b}
\bysame, \emph{Sur la stabilit\'e par produit tensoriel de complexes de
  {$\mathcal{D}$}-modules arithm\'etiques}, Manuscripta Math. \textbf{147}
  (2015), no.~1-2, 1--41.

\bibitem[Cha41]{Cha41}
C.~Chabauty, \emph{Sur les points rationnels des courbes alg{\'e}briques de
  genre sup{\'e}rieur {\`a} l'unit{\'e}}, C. R. Acad. Sci. Paris \textbf{212}
  (1941), 882--885.

\bibitem[Chi98]{Chi98}
B.~Chiarellotto, \emph{Weights in rigid cohomology applications to unipotent
  {F}-isocrystals}, Ann. Sci. Ecole. Norm. Sup. \textbf{31} (1998), 683--715.

\bibitem[CLS99a]{CLS99a}
B.~Chiarellotto and B.~Le~Stum, \emph{F-isocristaux unipotents}, Comp. Math.
  \textbf{116} (1999), 81--110.

\bibitem[CLS99b]{CLS99b}
\bysame, \emph{Pentes en cohomologie rigide et {$F$}-isocristaux unipotents},
  Manuscripta Mathematica \textbf{100} (1999), 455--468.

\bibitem[CT12]{CT12}
D.~Caro and N.~Tsuzuki, \emph{Overholonomicity of overconvergent
  ${F}$-isocrystals over smooth varieties}, Ann. of Math. \textbf{176} (2012),
  no.~2, 747--813.

\bibitem[Del70]{Del70}
P.~Deligne, \emph{Equations diff{\'e}rentielles {\`a} points singuliers
  r{\'e}guliers}, Lecture Notes in Mathematics, vol. 168, Springer-Verlag, New
  York, 1970.

\bibitem[Del89]{Del89}
\bysame, \emph{Le groupe fondamental de la droite projective moins trois
  points}, Galois Groups over $\mathbb{Q}$, Math. Sci. Res. Inst. Pub., 1989,
  pp.~79--297.

\bibitem[Del90]{Del90}
\bysame, \emph{Cat\'{e}gories tannakiennes}, The Grothendieck Festschrift, Vol.
  II, Progress in Mathematics, vol.~87, Birkh\"{a}user, 1990.

\bibitem[DMSS00]{DMSS00}
A.~Dimca, F.~Maaref, C.~Sabbah, and M.~Saito, \emph{Dwork cohomology and
  algebraic {D}-modules}, Math. Ann. \textbf{318} (2000), no.~1, 107--125.

\bibitem[EHS07]{EHS07}
H.~Esnault, P.~H. Hai, and X.~Sun, \emph{On {N}ori's fundamental group scheme},
  Geometry and Dynamics of Groups and Spaces, Progress in Mathematics, vol.
  265, Birkh{\"a}user, 2007, pp.~377--398.

\bibitem[Har75]{Har75}
R.~Hartshorne, \emph{On the de {R}ham cohomology of algebraic varieties}, Publ.
  Math. I.H.E.S. \textbf{45} (1975), no.~1, 6--99.

\bibitem[HJ10]{Had10}
M.~Hadian-Jazi, \emph{Motivic fundamental groups and integral points}, Ph.D.
  thesis, Universit{\"a}t Bonn, 2010.

\bibitem[HZ87]{HZ87}
R.~Hain and S.~Zucker, \emph{Unipotent variations of mixed {H}odge structure},
  Inventionnes Mathematica \textbf{88} (1987), 83--124.

\bibitem[Kim05]{Kim05}
M.~Kim, \emph{The motivic fundamental group of {$\P^1\setminus\{0,1,\infty\}$}
  and the theorem of {S}iegel}, Inventionnes Mathematica \textbf{161} (2005),
  629--656.

\bibitem[Kim09]{Kim09}
\bysame, \emph{The unipotent albanese map and {S}elmer varieties for curves},
  Publ. RIMS, Kyoto Univ. \textbf{45} (2009), 89--133.

\bibitem[Mac71]{Mac71}
S.~MacLane, \emph{Categories for the working mathematician}, Graduate Texts in
  Mathematics, vol.~5, Springer, 1971.

\bibitem[MD81]{MD81}
J.~Milne and P.~Deligne, \emph{Tannakian categories}, Hodge Cycles, Motives and
  Shimura Varieties, Lecture Notes in Mathematics, vol. 900, Springer, 1981,
  pp.~101--228.

\bibitem[NA93]{NA93}
V.~Navarro-Aznar, \emph{Sur la connection de {G}auss--{M}anin en homotopie
  rationelle}, Ann. Sci. Ecole. Norm. Sup. \textbf{26} (1993), 99--148.

\bibitem[Tam04]{Tam04}
A.~Tamagawa, \emph{Finiteness of isomorphism classes of curves in positive
  characterisitc with prescribed fundamental groups}, J. Algebraic Geom.
  \textbf{13} (2004), 675--724.

\bibitem[Wil97]{Wil97}
J.~Wildeshaus, \emph{Realisations of polylogarithms}, Lecture Notes in
  Mathematics, vol. 1650, Springer, 1997.

\end{thebibliography}

\end{document}